\documentclass[12pt]{article}
\usepackage[T1]{fontenc}
\usepackage[utf8]{inputenc}
\usepackage[english]{babel}
\usepackage{amssymb,amsmath}
\usepackage{amsthm}
\usepackage{cancel}
\usepackage{palatino}
\usepackage{color}
\usepackage{mathdots}
\usepackage{enumerate}
\numberwithin{equation}{section}
\newtheorem{theorem}{Theorem}[]     
\newtheorem{definition}[theorem]{Definition}
\newtheorem{proposition}[]{Proposition}
\newtheorem{lemma}[]{Lemma}
\newtheorem*{corollary}{Corollary}

\theoremstyle{definition}
\newtheorem*{example}{Example}
\newtheorem*{remark}{Remark}
\def\d{\partial}

\def\tb{\tilde{b}}

\def\tga{\tilde{\Gamma}}
\def\tna{\tilde{\nabla}}

\def\proof{\begin{center} {\bf Proof:} \end{center}\vspace{0.5pt}}
\def\QEDclosed{\mbox{\rule[0pt]{1.3ex}{1.3ex}}} 
\def\QED{\QEDclosed} 
\def\endproof{\hspace*{\fill}~\QED\par\endtrivlist\unskip}

\begin{document}

\title{Hamiltonian operators of Dubrovin-Novikov type in 2D}
\author{E.V. Ferapontov ${}^{*}$, P. Lorenzoni ${}^{**}$ and A. Savoldi ${}^{*}$}
    \date{}
    \maketitle
    \vspace{-7mm}
\begin{center}
${}^{*}$ Department of Mathematical Sciences, Loughborough University \\
Leicestershire LE11 3TU, Loughborough, United Kingdom \\
${}^{**}$ Dipartimento di Matematica e Applicazioni, University of Milano-Bicocca\\
via Roberto Cozzi 53 I-20125 Milano, Italy\\
e-mails: \\[1ex] \texttt{E.V.Ferapontov@lboro.ac.uk}\\
\texttt{paolo.lorenzoni@unimib.it}\\
\texttt{A.Savoldi@lboro.ac.uk}\\
\end{center}

\bigskip

\begin{abstract}

First order  Hamiltonian operators of differential-geometric type  were introduced by Dubrovin and Novikov in 1983, and thoroughly investigated by Mokhov. In 2D, they are generated by a pair of compatible flat metrics $g$ and $\tilde g$ which satisfy a set of additional  constraints coming from the skew-symmetry condition and the Jacobi identity.  We demonstrate that these constraints are equivalent to the requirement that $\tilde g$ is a linear Killing tensor of $g$ with zero Nijenhuis torsion. This allowed us to obtain a complete classification of $n$-component  operators with $n\leq 4$ (for $n=1, 2$ this was done before). 
For 2D operators the Darboux theorem does not hold: the operator may not  be  reducible to constant coefficient form. All interesting (non-constant) examples correspond to the case when the flat pencil $g, \tilde g$ is not semisimple, that is, the affinor $\tilde g g^{-1}$ has non-trivial Jordan block structure. In the  case of a direct sum of Jordan blocks with distinct eigenvalues we obtain a complete classification of Hamiltonian operators for any number of components $n$, revealing a remarkable correspondence with the class of trivial Frobenius manifolds modelled on $H^*({\bf CP}^{n-1})$.

\bigskip


\bigskip

\end{abstract}

\newpage

\tableofcontents

\newpage

\section{Introduction}

In 1983 Dubrovin and Novikov introduced Hamiltonian operators of differential-geometric type \cite{DN1}, 
\begin{equation}
P^{ij}= g^{ij}({\bf u}) \frac{d}{dx}+b^{ij}_k({\bf u})u^k_x,
\label{P}
\end{equation}
here ${\bf u}=(u^1, \dots, u^n)$ are the dependent variables,  and $i, j, k=1, \dots, n$ (we will assume   the non-degeneracy condition $\det g \ne 0$). The main observation was that the coefficients of such operators can be interpreted as differential-geometric objects: setting  $b^{ij}_k=-g^{is}\Gamma^j_{sk}$  and considering point transformations of the dependent variables, one can see that the coefficients $g^{ij}$ and $\Gamma^j_{sk}$ transform as components of a bivector (contravariant metric), and Christoffel's symbols of an affine connection, respectively. Imposing the requirement that the corresponding Poisson bracket, 
$$
\{F, G\}=\int \frac{\delta F}{\delta u^i} P^{ij}  \frac{\delta F}{\delta u^i}\  dx,
$$
is skew-symmetric and satisfies the Jacobi identity, one obtains that    the  bivector $g^{ij}$ defines a flat metric, and $\Gamma^j_{sk}$ is the associated Levi-Civita connection. This immediately establishes Darboux's theorem for such operators: in the flat coordinates of $g$ the operator $P$ takes constant coefficient form. Hamiltonian systems of hydrodynamic type are generated by Hamiltonians of the form $H=\int h({\bf u})dx$:
$$
u^i_t=P^{ij}\frac{\delta H}{\delta u^j}=\nabla^i \nabla_j h\ u^j_x.
$$
Such systems appear in a wide range of applications in hydrodynamics, chemical kinetics, the Whitham averaging method, the theory of Frobenius manifolds and so on, see the review papers \cite{DN, Tsarev} for further details and  references.

Hamiltonian operators of the form (\ref{P}) have  subsequently been generalised in a whole variety of different ways (degenerate, non-homogeneous, higher order,  multi-dimensional and non-local, see  \cite{Mokhov4} for a review), however, until now very few classification results are available due to the complexity of the problem. In this paper we 
address the  classification of  2D   Hamiltonian operators of  Dubrovin-Novikov type,
\begin{equation}
P^{ij}=g^{ij}({\bf u}) \frac{d}{dx}+b^{ij}_k({\bf u})u^k_x+\tilde g^{ij}({\bf u}) \frac{d}{dy}+
\tilde b^{ij}_k({\bf u})u^k_y,
\label{Ham}
\end{equation}
which are generated by a pair of flat metrics $g, \tilde g$, see \cite{Dubrovin1, Mokhov1, Mokhov2}. 
The operator $P$ will be called non-degenerate if the tensor $g+\lambda \tilde g$ is non-degenerate for generic  values of the parameter $\lambda$ (without any loss of generality we will assume  both $g$ and $\tilde g$ to be non-degenerate). Setting $b^{ij}_k=-g^{is}\Gamma^j_{sk}$ and $ \tilde b^{ij}_k=-\tilde g^{is}\tilde \Gamma^j_{sk}$,
where $\Gamma$ and $\tilde \Gamma$ are the Levi-Civita connections of  $g$ and $ \tilde g$, 
one introduces the obstruction tensor, 
$$
T^i_{jk}=\tilde \Gamma^i_{jk}- \Gamma^i_{jk}.
$$
 It is known that the vanishing of the obstruction tensor is necessary and sufficient for the existence of coordinates where the operator (\ref{Ham}) takes  constant coefficient form \cite{Dubrovin1} (2D Darboux theorem). Our analysis will be based on the following result of Mokhov:

\begin{theorem}\cite{Mokhov1} Let $g$ and $\tilde{g}$ be two flat metrics. Formula (\ref{Ham}) defines a Hamiltonian operator if and only if the obstruction tensor satisfies the relations
\begin{equation}
T^{ijk}=T^{kji},
\label{T1}
\end{equation}
\begin{equation}
T^{\{ijk\}}=0,
\label{T2}
\end{equation}
\begin{equation}
T^{ijs}T^r_{st}=T^{irs}T^j_{st},
\label{T3}
\end{equation}\begin{equation}
\nabla T^{ijk}=0,
\label{T4}
\end{equation}\begin{equation}
\tilde \nabla T^{ijk}=0.
\label{T5}
\end{equation}
Here $T^{ijk}=g^{ir}\tilde g^{ks}T^j_{rs}$, brackets $\{ \ \}$ denote cyclic permutations of the indices $i, j, k$, and $\nabla, \tilde \nabla$ are covariant differentiations in the Levi-Civita connections of $g, \tilde g$.
\end{theorem}
These relations imply that, in the flat coordinates of $g$, the second metric $\tilde g$  becomes linear, so that  the  classification of such operators reduces to the classification of algebras of certain type  \cite{Dubrovin1, Mokhov1}. This problem was addressed in \cite{Mokhov1}, resulting in a complete description of one- and two-component operators of the form (\ref{Ham}). Here we adopt a differential-geometric point of view: to proceed with the analysis of the above relations we introduce the affinor (that is, $(1,1)$-tensor) $L=\tilde g g^{-1}$ or, using indices, $L^i_j=\tilde g^{ik}g_{kj}$. According to \cite{Mokhov2},  $L$ must have zero Nijenhuis torsion  (the vanishing of the Nijenhuis torsion is a necessary condition for   1D brackets defined by $g$ and $\tilde g$ to be compatible \cite{Mokhov3, Fer}).

In the case when  $L$ has simple spectrum,
the results of \cite{Mokhov2} imply the existence of coordinates where the Hamiltonian operator $P$ takes constant coefficient form. It turns out that all interesting (non-constant) examples correspond to the case when $L$ has non-trivial Jordan block structure. The simplest known example of this kind is provided by the two-component operator
\begin{equation}\label{2d_operator}
P=
\begin{pmatrix}
0 & 1 \\
1 & 0
\end{pmatrix}
\frac{d}{dx}
+\begin{pmatrix}
-2u^1 & u^2 \\
u^2 & 0
\end{pmatrix}
\frac{d}{dy}
+\begin{pmatrix}
-u_y^1 & 2 u_y^2 \\
-u_y^2 & 0
\end{pmatrix},
\end{equation}
which is related to the Lie algebra of vector fields on the plane \cite{Dubrovin1, Mokhov1}. It is generated by the flat contravariant metrics
$$
g=\begin{pmatrix} 
0 & 1 \\
1 & 0
\end{pmatrix}
, ~~~ \tilde g=\begin{pmatrix}
-2u^1 & u^2 \\
u^2 & 0
\end{pmatrix}.
$$
One can easily see that, for generic values of $u^1, u^2$, the corresponding affinor $L=\tilde g g^{-1}$ is a single $2\times 2$ Jordan block. 



\section{Summary of the main results}

Our first result establishes a link between 2D Hamiltonian operators (\ref{Ham}) and the theory of Killing tensors:

\begin{theorem}\label{Kill_thm} Let $g$ and $\tilde{g}$ be two flat metrics which define the Hamiltonian operator (\ref{Ham}). The Mokhov conditions \eqref{T1}--\eqref{T5} are equivalent to the following:
\begin{enumerate}
\item Linearity of the bivector $\tilde{g}^{ij}$ in the flat coordinates of $g$. Invariantly, this means $\nabla^2 \tilde g=0$ where $\nabla$ denotes covariant differentiation in the Levi-Civita connection of $g$ (this fact was established earlier in \cite{Dubrovin1, Mokhov1}).
\item The vanishing of the Nijenhuis torsion of the affinor $L^i_j=\tilde{g}^{il} g_{lj}$ \cite{Mokhov2}.
\item The Killing condition for the bivector $\tilde g$: 
$
\nabla^i\tilde{g}^{kj}+\nabla^k\tilde{g}^{ij}+\nabla^j\tilde{g}^{ik}=0.
$
\end{enumerate}
In particular, the flatness of $g$ and the above three conditions imply the flatness of the second metric $\tilde{g}$.
\end{theorem}
Thus, the classification of Hamiltonian operators of the form (\ref{Ham}) is reduced to the classification of linear Killing bivectors with zero Nijenhuis torsion in flat pseudo-Euclidean spaces. A  tensorial proof of Theorem 2 is given in Section \ref{sect_killing}. Using the fact that any  Killing bivector in  flat space is the sum of symmetrised tensor products of Killing vectors, we obtain a complete classification of 2D Hamiltonian operators with $n\leq 4$ components (Section 5).

The Killing condition  plays also a key role in the proof of the splitting property for Hamiltonian operators, which can be seen as an analogue of the splitting lemma for affinors with zero Nijenhuis torsion proved by Bolsinov and Matveev \cite{BM} in the context of projectively equivalent metrics. First of all we recall their result. Let $M$ be an $n$-dimensional manifold,  and let $L$ be an affinor on $M$ with zero Nijenhuis torsion. Suppose that there exists a frame (not necessarily holonomic) in which  $L$ takes block diagonal form,
\begin{equation}\label{affinor_split}
L=
\begin{pmatrix}
A & 0\\
0 & B
\end{pmatrix},
\end{equation}
where $\mathrm{Spec}(A) \, \cap \, \mathrm{Spec}(B) = \emptyset$. Then there exists a coordinate system $({\bf u},{\bf v})=(u^1,\ldots,u^m, v^{m+1},\ldots, v^n)$ such that $A$ depends on ${\bf u}$ and $B$ depends on ${\bf v}$ only, that is, $L$ is a direct sum of two affinors (both with vanishing Nijenhuis torsion). Adding the Killing condition, in Section \ref{sect_splitting} we  show how to extend this splitting structure to the metrics, namely we  prove that, in the same coordinate system,  the two metrics $g$ and $\tilde g$ also assume block diagonal forms,
$$
g=
\begin{pmatrix}
g_1({\bf u}) & 0\\
0 & g_2({\bf v})
\end{pmatrix},
\quad \quad
\tilde g=
\begin{pmatrix}
\tilde{g_1}({\bf u}) & 0\\
0 & \tilde{g_2}({\bf v}) 
\end{pmatrix}.
$$
This suggests the definition of reducible operators:
given an $m$-component operator $P_1$ with the dependent variables $u^1, \dots, u^m$, and an $(n-m)$-component operator $P_2$ with the dependent variables  $v^{m+1}, \dots, v^{n}$,  their direct sum is  the $n$-component operator $P$ defined by the formula 
$$
P=\left(\begin{array}{cc}
P_1&0\\
0&P_2
\end{array}
\right),
$$
on the combined set of variables $(u^1, \ldots, u^m, v^{m+1}, \ldots, v^n)$. The corresponding metrics $g$, $\tilde g$ will be  direct sums of the metrics defining $P_1$ and $P_2$.
Operators of this type will be called \emph{reducible}. Thus, our second result can be formulated as follows:

\medskip

\noindent {\bf The Splitting Lemma.\ }
{\it Let $P$ be a Hamiltonian operator such that the corresponding affinor $L=\tilde g g^{-1}$ can be represented in the block-diagonal form (\ref{affinor_split}) in some (non-holonomic) frame, and let $\mathrm{Spec}(A) \, \cap \, \mathrm{Spec}(B) = \emptyset$. Then $P$ decouples into a direct sum of two Hamiltonian operators, with the corresponding affinors $A$ and $B$.}

\medskip

Thus, any Hamiltonian operator (\ref{Ham}) can be represented as a direct sum of irreducible   operators $P_{\alpha}$
(each generated by a pair of flat metrics $g_{\alpha}$, $\tilde g_{\alpha},$ defined on a manifold of dimension $n_{\alpha}$) such that the corresponding affinor $L_{\alpha}=\tilde g_{\alpha}g_{\alpha}^{-1}$ either has a unique real eigenvalue of  multiplicity $n_{\alpha}$, or a pair of complex conjugate eigenvalues of the same multiplicity (in the last case $n_{\alpha}$ must be even).

As a consequence of the splitting lemma we will prove that,  if 
the affinor $L$ is diagonal, then
the Hamiltonian operator can be brought to   constant coefficient form. This generalises the analogous result of  \cite{Mokhov2} obtained under the additional assumption of  the simplicity of the spectrum of $L$. In what follows, we will be interested in Hamiltonian operators which are not reducible, and not transformable to  constant coefficient form.

Our approach to the classification of Hamiltonian operators in 2D is based on the Killing property. As an illustration, in  Section \ref{classification} we review the already known two-component case \cite{Mokhov2}, and   give a complete classification of three- and four-component Hamiltonian operators in 2D. In the three-component case, the  main result is as follows.
\begin{theorem}\label{thm3}
Any  irreducible non-constant three-component Hamiltonian operator in 2D can be brought (by a change of the dependent variables $u^i$) to the form $\pm P$ where $P$ can have one of the following two canonical forms (in both cases the  affinor $L$ is a single $3\times 3$ Jordan block):
\begin{enumerate}
\item
Jordan block with constant eigenvalue
$$
P=
\begin{pmatrix}
0 & 0 & 1 \\
0 & 1 & 0 \\
1 & 0 & 0
\end{pmatrix}
\frac{d}{dx}
+\begin{pmatrix}
-2 u^2 & u^3 & \lambda \\
u^3 & \lambda & 0 \\
\lambda & 0 & 0
\end{pmatrix}
\frac{d}{dy}
+\begin{pmatrix}
-u_y^2 & 2 u_y^3 & 0\\
-u_y^3 & 0 & 0\\
0 & 0 & 0
\end{pmatrix},
$$
\item
Jordan block with non-constant eigenvalue
$$
P=
\begin{pmatrix}
0 & 0 & 1 \\
0 & 1 & 0 \\
1 & 0 & 0
\end{pmatrix}
\frac{d}{dx}
+\begin{pmatrix}
-2u^1 & -\frac{1}{2}u^2 & u^3 \\
-\frac{1}{2}u^2 & u^3 & 0 \\
u^3 & 0 & 0
\end{pmatrix}
\frac{d}{dy}
+\begin{pmatrix}
-u_y^1 &\frac{1}{2} u_y^2 & 2u_y^3\\
-u_y^2 & \frac{1}{2} u_y^3 & 0\\
-u_y^3 & 0 & 0
\end{pmatrix}.
$$
\end{enumerate}
\end{theorem}
\begin{remark}
The second example is a three-component version of the $n$-component Hamiltonian operator introduced by Mokhov in \cite{Mokhov1}.
\end{remark}

In the four-component  situation calculations become more complicated, and we get several canonical forms labelled by Segre types of the affinor $L$, see end of Section \ref{classification}.

Although  our approach works for any number of components $n$, for $n>4$  computations become rather cumbersome. The main difficulty is  when the affinor $L$ consists of several Jordan blocks with the same eigenvalue. In Section \ref{sect_Jordan}
we analyse the particular case of a single $n\times n$ Jordan block: 
\begin{theorem}\label{Jor_thm}
Let $P$ be a Hamiltonian operator (\ref{Ham}) such that the affinor  $L= \tilde{g} g^{-1}$ is a single $n\times n$ Jordan block with non-constant eigenvalue. Then there exists a coordinate system in which $g$  and $\tilde g$ can be reduced to the following canonical forms:
$$
g=\pm
\begin{pmatrix}
 &&  1\\
& \iddots &  \\
1 & & 
\end{pmatrix},
\quad
\tilde{g}=\pm
\left\{
\begin{array}{lcl}
\mu^{(n;0)} & \mbox{if} & n\not\equiv 1\!\!\!\! \mod 3,\\
\mu^{(n;0)} + \kappa_1 \mu^{(n;\frac{n-1}{3})} & \mbox{if} & n\equiv 1\!\!\!\! \mod 3, n \neq 4,\\
\mu^{(4;0)}+\kappa_1 \mu^{(4;1)}+{\mu} & \mbox{if} & n=4.
\end{array}
\right.
$$
Here $ \kappa_1$ is an arbitrary constant,  the symmetric bivector $\mu^{(n;k)}$ is defined as
$$
\mu^{(n;k)ij}=[3(i+j)-2(n+2-k)]u^{i+j-1+k},
$$
and $\mu$ is the constant symmetric matrix
${\mu}^{ij}=\delta^{i,4-j}+\lambda \delta^{i,5-j}$, $\lambda=const.$

\end{theorem}
For the constant eigenvalue case, the  analogous statement  is given at the end of  Section 6. The bivector $\mu^{(n;0)}$  corresponds to the Mokhov operator \cite{Mokhov1}, which is of the single Jordan block type for all $n\neq 4$.  Thus, for $n=3$, the bivector $\mu^{(3;0)}$  gives rise to  the second case of Theorem 3. Note that, for $n= 4$,  the Mokhov operator has the type of two $2\times 2$ Jordan blocks: this explains the presence of the term $\mu$ in the formula for $\tilde g$ in Theorem 4. This result, combined with the splitting lemma,
 provides a complete classification of 2D operators of Dubrovin-Novikov type in the case of a direct sum of Jordan blocks with distinct eigenvalues.  

In Section \ref{sect_frobenius} we  show that the case of a single $n\times n$ Jordan block  with non-constant eigenvalue gives rise to the \emph{trivial} non-semisimple Frobenius manifold  whose underlying Frobenius algebra corresponds to the cohomology ring of ${\bf CP}^{n-1}$.

Finally, in Section 7 we  extend our approach to  Hamiltonian operators in dimensions higher than two. Recall that the results of  \cite{Mokhov1} imply that, for $d\geq 3$, there exist no non-trivial $d$-dimensional one- or two-component Hamiltonian operators of the form
$$
P^{ij}=\sum_{\alpha =1}^d \left(g^{ij\alpha}({\bf u}) \frac{d}{dx^{\alpha}}+b^{ij\alpha}_k({\bf u})u^k_{x^{\alpha}}\right),
$$
in other words, any one- or two-component operator of this kind can be transformed to constant coefficient form. 
We obtain a complete description of three-component operators which are essentially three-dimensional, and cannot be transformed to constant coefficients:
\begin{theorem} Any non-degenerate three-component Hamiltonian operator in 3D, which is not transformable to constant coefficients,  can be brought to one of the two canonical forms:
$$
P=
\begin{pmatrix}
\partial_z & 0 & \partial_x \\
0 & \partial_x & 0 \\
\partial_x & 0 & 0
\end{pmatrix}
+
\begin{pmatrix}
-2u^2 \partial_y -u_y^2 & u^3 \partial_y +2 u_y^3 & 0 \\
u^3 \partial_y -u_y^3 & 0 & 0\\
0 & 0 & 0
\end{pmatrix},
$$
or
$$
P=
\begin{pmatrix}
0 & \partial_x & 0 \\
\partial_x & 0 & 0 \\
0 & 0 & \partial_z
\end{pmatrix}
+
\begin{pmatrix}
-2 u^1 \partial_y -u_y^1  & u^2\partial_y +2 u_y^2 & 0 \\
u^2 \partial_y -u_y^2& 0 & 0\\
0 & 0 & 0
\end{pmatrix}.
$$
(here we allow arbitrary changes of the dependent variables $u^i$, and linear transformations of the independent variables $x, y, z$).
\end{theorem}
Note that the second operator is reducible: it is a direct sum of the non-constant two-component operator (\ref{2d_operator}), and the operator $\partial_z$. Thus, there exists a unique irreducible three-component operator in 3D.


\section{Hamiltonian operators  and linear Killing tensors with zero Nijenhuis torsion: proof of Theorem \ref{Kill_thm}}\label{sect_killing}
In this section we  rewrite Mokhov's conditions \eqref{T1} -- \eqref{T5} in the form which is more suitable for our purposes, making link with the theory of Killing tensors. 


\vspace{5pt}\noindent
{\bf Theorem \ref{Kill_thm}.} {\it
Let $g$ be a flat metric.  Then conditions  \eqref{T1} -- \eqref{T5} are equivalent to the following:
\begin{enumerate}
\item Linearity of $\tilde{g}^{jk}$ in the flat coordinates of $g$.
\item Vanishing of the Nijenhuis torsion of the affinor $L^i_j=\tilde{g}^{il} g_{lj}$.
\item The Killing condition: 
\begin{equation}\label{killing}
\nabla^i\tilde{g}^{kj}+\nabla^k\tilde{g}^{ij}+\nabla^j\tilde{g}^{ik}=0.
\end{equation}
\end{enumerate}
In particular, the flatness of $g$ and the above three conditions imply the flatness of the second metric $\tilde{g}$.
}

\begin{remark}
The facts  that $\tilde{g}$ must be linear in the flat coordinates of $g$, and that the Nijenhuis torsion of $L$ must vanish,
are well known  \cite{Dubrovin1, Mokhov1, Mokhov2}. They are equivalent to \eqref{T4} and \eqref{T1}, respectively. Our contribution here is the Killing property, and the observation that the assumption of  flatness of $\tilde g$ can be dropped.
\end{remark}

\begin{center} {\bf Proof of Theorem \ref{Kill_thm}:} \end{center}\vspace{0.5pt}

\noindent
{\bf(a)}. \emph{The condition \eqref{T1} is equivalent to the vanishing of the Nijenuis torsion of $L$}.

\noindent
This  was proved by Mokhov \cite{Mokhov3,Mokhov5}, here we briefly recall the proof.
Let $\tilde{b}^{ij}_k=-\tilde{g}^{is}\tilde{\Gamma}^{j}_{sk}$ be  contravariant Christoffel symbols of the second metric,
 by definition they satisfy the 
 conditions
\begin{eqnarray*}
\d_k\tilde{g}^{ij}&=&\tilde{b}^{ij}_k+\tilde{b}^{ji}_k,\\
\tilde{g}^{il}\tilde{b}^{jk}_l&=&\tilde{g}^{jl}\tilde{b}^{ik}_l.
\end{eqnarray*}
 Written in the flat coordinates of $g$, the condition \eqref{T1} reads
$$g^{il}\tilde{b}^{jk}_l=g^{jl}\tilde{b}^{ik}_l.$$
Thus, the metrics  $\tilde{g}$ and $g$ are almost compatible, and this is known to be equivalent to the vanishing of the Nijenhuis torsion  \cite{Mokhov3}.


\noindent
{\bf(b)}. \emph{The condition \eqref{T2} is equivalent to the Killing property}.

\noindent
Using \eqref{T1} we can rewrite \eqref{T2} as
\begin{eqnarray*}
&&\sum_{(i,j,k)}[T^{ijk}+T^{kji}]=0.
\end{eqnarray*}
In the flat coordinates of $g$ we have
\begin{eqnarray*}
&&\sum_{(i,j,k)}[T^{ijk}+T^{kji}]=\\
&&=\tilde{g}^{ks}g^{ir}\tilde{\Gamma}^{j}_{rs}+\tilde{g}^{is}g^{kr}\tilde{\Gamma}^{j}_{rs}+
\tilde{g}^{is}g^{jr}\tilde{\Gamma}^{k}_{rs}+\tilde{g}^{js}g^{ir}\tilde{\Gamma}^{k}_{rs}+
\tilde{g}^{js}g^{kr}\tilde{\Gamma}^{i}_{rs}+\tilde{g}^{ks}g^{jr}\tilde{\Gamma}^{i}_{rs}=\\
&&=-[g^{ir}\tilde{b}^{kj}_{r}+g^{kr}\tilde{b}^{ij}_{r}+
g^{jr}\tilde{b}^{ik}_{rs}+g^{ir}\tilde{b}^{jk}_{s}+
g^{kr}\tilde{b}^{ji}_{s}+g^{jr}\tilde{b}^{ki}_{r}]=\\
&&=-[g^{ir}\d_r\tilde{g}^{kj}+g^{kr}\d_r\tilde{g}^{ij}+
g^{jr}\d_r\tilde{g}^{ik}]=\\
&&=-[\d^i\tilde{g}^{kj}+\d^k\tilde{g}^{ij}+\d^j\tilde{g}^{ik}]=0.
\end{eqnarray*}
In invariant notation, this gives  the \emph{Killing condition},
$$\nabla^i\tilde{g}^{kj}+\nabla^k\tilde{g}^{ij}+\nabla^j\tilde{g}^{ik}=0,$$
here $\nabla$ is the Levi-Civita connection of $g$.


\noindent
{\bf(c)}. \emph{The condition \eqref{T4} is equivalent to the linearity of $\tilde{g}$ in the flat coordinates of $g$}.

\noindent
 In the flat coordinates
 of $g$,  \eqref{T4} implies
\begin{eqnarray*} 
\d_r(T^{ijk}+T^{ikj})=\d_r[g^{it}(\tilde{b}^{kj}_{t}+\tilde{b}^{jk}_{t})]
=\d_r\d^i\tilde{g}^{jk}=0.
\end{eqnarray*}
This means that $\tilde{g}$ is linear. 
Conversely,  assuming that $\tilde{g}$ is linear in the flat coordinates of $g$, and using  \eqref{T1} and \eqref{T2}, we obtain \eqref{T4}:
$$
0=\d_r(T^{ijk}+T^{ikj})=\d_r(T^{ijk}+T^{jki})=-\d_r T^{kij}.
$$


\noindent
{\bf(d)}. \emph{The conditions \eqref{T3} and \eqref{T4} are equivalent to the flatness of $\tilde{g}$}.

\noindent
The condition \eqref{T4} means that, in the flat coordinates of $g$, the contravariant Christoffel symbols
 $\tilde{b}^{ij}_k$ are constant. This follows from the identity
 $$-\d_rT^{kij}=\d_r(g^{km} \tilde{b}^{ij}_m)=g^{km} \d_r\tilde{b}^{ij}_m=0.$$
 Similarly,   the condition \eqref{T5} means that, in the flat coordinates of $\tilde{g}$, the contravariant Christoffel symbols
 $b^{ij}_k$ are constant. 
 Written in the flat coordinates of $g$,  the condition \eqref{T3} reads
$$
\tilde{g}^{sq}g^{ip}\tilde{\Gamma}^{j}_{pq}\tilde{\Gamma}^{r}_{st}=
\tilde{g}^{sq}g^{ip}\tilde{\Gamma}^{r}_{pq}\tilde{\Gamma}^{j}_{st},
$$
which is equivalent to
$$\tilde{b}^{sj}_{p}\tilde{b}^{ir}_{s}-\tilde{b}^{sr}_{p}\tilde{b}^{ij}_{s}=0.$$
Due to  \eqref{T4},  the vanishing of the curvature of $\tilde{\nabla}$, written in the flat coordinates of $g$,  reads
\begin{eqnarray*}
g^{is}\left(\partial_s\tilde{b}^{jr}_p-\partial_p\tilde{b}^{jr}_s\right)
-\tilde{b}^{ij}_s\tilde{b}^{sr}_p+\tilde{b}^{ir}_s\tilde{b}^{sj}_p= 
-\tilde{b}^{ij}_s\tilde{b}^{sr}_p+\tilde{b}^{ir}_s\tilde{b}^{sj}_p=0.
\end{eqnarray*}

\noindent
{\bf(e)}. \emph{The condition \eqref{T5} can be dropped}.

\noindent
We recall that, in the flat coordinates of $g$, we have $T^{lkj}=-g^{lm} \tb^{kj}_m$ and $T^i_{jk}=\tga^i_{jk}=\tga^i_{kj}=T^i_{kj}$ (by the symmetry of  $\tilde \nabla$). A straightforward computation gives
\begin{eqnarray*}
\tna_r T^{ijk}=-(T^i_{rl} T^{lkj}+T^k_{rl} T^{ilj} +T^j_{rl} T^{ikl}).
\end{eqnarray*}
Using conditions \eqref{T1}, \eqref{T2} and \eqref{T3} we  obtain
\begin{eqnarray*}
\tna_r T^{ijk}=-T^k_{rl} (-T^{lji} +T^{ijl}),
\end{eqnarray*}
where the last term vanishes by \eqref{T1}.

\noindent
{\bf(f)}. \emph{The flatness of $\tilde g$ follows from the  flatness of $g$,  linearity of $\tilde{g}$, the Killing condition, and the vanishing of the Nijenhuis torsion}.

\noindent
Since $\tilde{g}^{ij}=(\tilde{b}^{ij}_l+\tilde{b}^{ji}_l)u^l+g_0^{ij}$, the Killing condition reads
$$
g^{is}(\underline{\tilde{b}^{kj}_s}+\underline{\tilde{b}^{jk}_s})+g^{ks}(\tilde{b}^{ij}_s+\underline{\tilde{b}^{ji}_s})
+g^{js}(\tilde{b}^{ik}_s+\tilde{b}^{ki}_s)=0.
$$ 
Then, using  \eqref{T1} for the underlined terms, we can rewrite the Killing condition as
\begin{equation}\label{fc}
\tilde{b}^{ij}_{s}g^{sk}+(\tilde{b}^{ki}_{s}+\tilde{b}^{ik}_{s})g^{sj}=0.
\end{equation}
We will make use of this condition later in Section 6.1.
In the flat coordinates of $g$, the vanishing of the Nijenhuis torsion $\mathcal{N}(L)$ of  the affinor $L$ reads
\begin{gather*}
0=\mathcal{N}^k_{ij}=L^s_i\d_sL^k_j-L^s_j\d_sL^k_i+L^k_s\d_jL^s_i-L^k_s\d_iL^s_j\\
=\tilde{g}^{sl}g_{li}g_{mj}(\tilde{b}^{km}_s+\tilde{b}^{mk}_s)-\tilde{g}^{sl}g_{lj}g_{mi}
(\tilde{b}^{km}_s+\tilde{b}^{mk}_s)
+\tilde{g}^{kl}g_{ls}g_{mi}(\tilde{b}^{sm}_j+\tilde{b}^{ms}_j)\\
-\tilde{g}^{kl}g_{ls}g_{mj}(\tilde{b}^{sm}_i+\tilde{b}^{ms}_i).
\end{gather*}
Multiplying by $g^{ip} g^{jq}$, taking the sum over $i$ and $j$ and using  \eqref{fc} we get
\begin{gather}
g^{ip} g^{jq}\mathcal{N}^k_{ij}
=\tilde{g}^{sp}(\tilde{b}^{kq}_s+\tilde{b}^{qk}_s)-\tilde{g}^{sq}(\tilde{b}^{kp}_s+\tilde{b}^{pk}_s)+\tilde{g}^{ks}(\tilde{b}^{qp}_s-\tilde{b}^{pq}_s)=:J^{pqk}.\label{Nij2}
\end{gather}
Thus, $\mathcal{N}(L)=0$ if and only if \eqref{Nij2} holds, that is, if and only if $J^{pqk}=0$. Let us now consider the sum $J^{pqk}+J^{qkp}$: it must be zero due to the vanishing of the Nijenhuis torsion. A direct computation gives
\begin{gather}
0=J^{pqk}+J^{kpq}
=2(\tilde{g}^{ks}\tilde{b}^{qp}_s-\tilde{g}^{qs}\tilde{b}^{kp}_s).\label{Nij3}
\end{gather}
Assume now that $\tilde{g}$ is linear in the flat coordinates of $g$,
\begin{equation}\label{lm}
\tilde{g}^{ij}=c^{ij}_k u^k+g_0^{ij}. 
\end{equation}
This implies
\begin{equation}\label{Beqb}
\partial_k \tilde{g}^{ij}=c^{ij}_k=\tilde{b}^{ij}_k +\tilde{b}^{ji}_k.
\end{equation}
For the linear metric \eqref{lm}, the condition \eqref{Nij3}  reads
$$(c^{ks}_l u^l+g_0^{ks})\tilde{b}^{qp}_s=(c^{qs}_l u^l+g_0^{qs})\tilde{b}^{kp}_s.$$
This is equivalent to
$$c^{ks}_l \tilde{b}^{qp}_s=c^{qs}_l\tilde{b}^{kp}_s
\quad \mbox{ and } \quad 
g_0^{ks}\tilde{b}^{qp}_s=g_0^{qs}\tilde{b}^{kp}_s.$$
Due to \eqref{Beqb} the first condition can be written as
$$(\tilde{b}^{ks}_l+\tilde{b}^{sk}_l)\tilde{b}^{qp}_s=(\tilde{b}^{qs}_l+\tilde{b}^{sq}_l)\tilde{b}^{kp}_s.$$
Using  \eqref{fc} we obtain
$$g_{m r} (\underline{g^{s j}}\tilde{b}^{k r}_{j}\underline{\tilde{b}^{qp}_s}-g^{qj}\tilde{b}^{s r}_{j}\tilde{b}^{kp}_s)=0,$$
and using  \eqref{T1} for the underlined terms, we finally get
$$g_{m r} g^{qs}(\tilde{b}^{j p}_s \tilde{b}^{k r}_{j}-\tilde{b}^{j r}_{s}\tilde{b}^{kp}_{j})=0.$$

\endproof

\begin{remark}
Using Mokhov's conditions it is easy to prove that $g$ and the homogeneous linear part of $\tilde{g}$
 define an exact flat pencil  of metrics. More precisely, we have
\begin{eqnarray*}
\mathcal{L}_{X} g^{ij}&=&g_1^{ij},\\
\mathcal{L}_{X} g_1^{ij}&=&0,
\end{eqnarray*}
where $g_1^{ij}=(\tilde{b}^{ij}_l+\tilde{b}^{ji}_l)u^l$ and $X^i=-g_1^{is}g_{sl}u^l$. Moreover, 
 $X$ is constant in flat coordinates of $g_1$ ($\nabla_1X=0$). Exactness is one of the main properties of flat  pencils of metrics related to Frobenius manifolds. In this case $X$
 is the vector field defining the unit of the multiplicative structure. This observation suggests that flat pencil of metrics defining
  2D Hamiltonian operators might be related to Frobenius manifolds. We will discuss this point in Section 6
   in the case of Mokhov's example.
\end{remark}

\begin{remark}
In what follows we will need an alternative form of the Killing condition, namely
\begin{equation}\label{kill1}
g^{is}\partial_s\tilde{g}^{kj}+g^{ks}\partial_s\tilde{g}^{ij}+g^{js}\partial_s\tilde{g}^{ik}-
\tilde{g}^{is}\partial_s g^{kj}-\tilde{g}^{ks}\partial_s g^{ij}-\tilde{g}^{js}\partial_s g^{ik}=0.
\end{equation}
This  can be obtained as follows.
Computing  covariant derivative of $\tilde g^{ij}$ we get
\begin{gather}\label{kill_contravariant}
g^{ks}\nabla_s\tilde{g}^{ij}
=g^{ks}\partial_s \tilde{g}^{ij} -b^{ki}_{m}\tilde{g}^{mj} -b^{kj}_{m}\tilde{g}^{im}.
\end{gather}
Using $
\partial_s g^{ij}= b^{ij}_s+b^{ji}_s
$  and substituting \eqref{kill_contravariant} into the Killing conditions, one arrives at \eqref{kill1}.
\end{remark}


\section{The splitting lemma}\label{sect_splitting}

One of the conditions which follows from the Hamiltonian property is the vanishing of the Nijenhuis torsion of  the affinor $L=\tilde g g^{-1}$.
In the hypothesis that the spectrum of  $L$ can be decoupled into two  subsets with  empty intersection,  Bolsinov and Matveev  established the following splitting property:

\begin{lemma}[Splitting Lemma, \cite{BM1}]
\label{lem4}
Let $L$ be an affinor  with zero Nijenhuis torsion on a manifold $M$, $\mathrm{dim} M =n$. Suppose  there exists a (non-holonomic) frame in which $L$ takes block diagonal form,   
$$
L=
\begin{pmatrix}
A& 0 \\
0 & B
\end{pmatrix},
$$
where $\mathrm{Spec}(A)\, \cap \, \mathrm{Spec}(B)=\emptyset$. 
Then there exists a local coordinate system $(u^1, ..., u^m,$ $ v^{m+1}, ..., v^n)$ such
that $$
L=\begin{pmatrix} A({\bf u}) & 0 \\ 0 & B({\bf v}) \end{pmatrix}.
$$
\end{lemma}

Using the Killing condition, one can extend the splitting structure  to the metrics.
First of all we recall a well-known fact from linear algebra: in the hypothesis of the above lemma, if $g$ and $\tilde{g}$ are two non-degenerate symmetric bivectors  related by the affinor $L$, that is $\tilde g^{ij}:=L^j_k g^{ki}$, then $g$ and $\tilde g$ assume the form
\begin{equation}\label{split1}
g=
\begin{pmatrix}
\sigma & 0\\
0 & \eta
\end{pmatrix},
\quad \quad
\tilde g=
\begin{pmatrix}
\tilde{\sigma} & 0\\
0 & \tilde{\eta} 
\end{pmatrix},
\end{equation}
(see \cite{BM1} for  applications of this result to the theory of projectively equivalent metrics).

\begin{lemma}\label{split}
In the hypothesis of Lemma \ref{lem4}, let $g$ and $\tilde g$ be two non-degenerate symmetric bivectors \eqref{split1} such that  $\tilde g^{ij}:=L^j_k g^{ki}$. 
If the Killing condition \eqref{killing} holds, then $\sigma$, $\tilde{\sigma}$ must depend only on ${\bf u} =(u^1, \ldots, u^m)$,  and $\eta$, $\tilde{\eta}$ must depend only on ${\bf v} =(v^{m+1}, \ldots, v^n)$.
\end{lemma}

\proof
By Lemma \ref{lem4},   $A=A({\bf u})$ is an $m\times m$ matrix,  and $B=B({\bf v})$ is an $(n-m) \times (n-m) $ matrix.
Let  $I=\{1, \ldots, m\}$ and $J=\{m+1, \ldots, n\}$. We know  that if $i \in I$ and $j \in J$, then $g^{ij}=0$.
Then, for $i\in I$ and $j,k\in J$, the condition \eqref{kill1} leads to
$$
g^{is}\partial_s\tilde{g}^{kj}-\tilde{g}^{is}\partial_s g^{kj}=0,
$$
 in particular,
$$
\sigma^{is}\partial_s\tilde{\eta}^{kj}-\tilde{\sigma}^{is}\partial_s \eta^{kj}=0.
$$
Multiplying by the inverse matrix $\sigma_{li}$ we obtain
$$
\partial_l (B^k_p \eta^{pj})-A^s_l \partial_s \eta^{kj}=0,
$$
as $\tilde{\sigma}^{si}\sigma_{il}=A^s_l$. Since $l \in I$ and the elements of $B$ depend  on $\bf {v} $ only, our relation becomes
$$
B^k_p  \partial_l \eta^{pj}-A^s_l \partial_s \eta^{kj}=0.
$$
Fixing $j$, let $C^{k}_i:=\partial_i \eta^{kj}$, then we get $B^k_p C^p_l=C^k_s A^s_l  $, that is
$$
B C=C A.
$$
As $\mathrm{Spec}(A)  \, \cap \,  \mathrm{Spec}(B)=\emptyset$, it follows  that $C \equiv 0$. 
Thus 
$$
\partial_i \eta^{jk}=0, \ \forall\; i \in I,  \ \forall\; j,k \in J.
$$
If we now take $i \in J$ and $j,k \in I$, following the same method we get
$$
\partial_i \sigma^{jk}=0, \ \forall\; i \in J,  \ \forall\; j,k \in I.
$$
\endproof

This establishes the Splitting Lemma formulated in Section 2. It allows us to focus on affinors with one single eigenvalue, otherwise we can split them and consider each block separately. 


As a simple application of the splitting lemma we can establish Darboux's theorem for Hamiltonian operators whose affinor $L$ is diagonal (has no non-trivial Jordan blocks: note that we allow coinciding eigenvalues). It is based on the following result:

\begin{proposition}
Let $L$ be a diagonal affinor, $g$ be a flat contravariant metric, and $\tilde g= L g$. Suppose that the Nijenhuis torsion of $L$ vanishes, and the Killing condition holds. Then there exists a coordinate system where $L$ and $g$ take constant coefficient form.
\end{proposition}
\proof
Since the Nijenhuis torsion of $L$ vanishes, using the splitting lemma we can bring $L$ to block diagonal form,
$$
L=
\begin{pmatrix}
L_1 & & &\\
& L_2 & &\\
& & \ddots & \\
& & & L_k
\end{pmatrix}.
$$
Here each $L_i$ is a scalar operator with the same eigenvalue,
$$
L_i=
\begin{pmatrix}
\lambda^i & &\\
& \ddots& \\
& & \lambda^i
\end{pmatrix},
$$
$\lambda^i\neq \lambda^j$ for $i\neq j$, and $\lambda^i$ depends  on coordinates of its own block only. By the splitting lemma for the metrics, we have
$$
g=
\begin{pmatrix}
g_{\lambda^1} & & \\
& \ddots& \\
& & g_{\lambda^k}
\end{pmatrix},
\quad
\tilde g=
\begin{pmatrix}
\lambda^1 g_{\lambda^1} & & \\
& \ddots& \\
& & \lambda^k g_{\lambda^k}
\end{pmatrix},
$$
where $g_{\lambda^i}$  depends  on coordinates of its own block only.
Thus we can consider  each block separately. For instance, suppose $L_1$ is an $m\times m$ scalar operator with the eigenvalue $\lambda_1$. Let us set $\lambda=\lambda^1$ and $h=g_{\lambda^1}$. We know  that $\lambda$ and  $h$ depend  on $u^1, \ldots u^m$ only,  and no other block depends on these coordinates.
The condition \eqref{kill1} leads to
$$
h^{kj} h^{is}\partial_s \lambda+h^{ji} h^{ks}\partial_s\lambda+h^{ik}h^{js}\partial_s\lambda=0.
$$
Since $h$ is non-degenerate, contracting with $h_{qi} h_{pj}$ we get
$$
\delta^k_p \partial_q \lambda+ h_{pq} h^{ks}\partial_s\lambda+\delta^{k}_q \partial_p\lambda=0.
$$
Setting $q=k$ and summing over $k$ we obtain
$$
 \partial_p \lambda+ h_{pk} h^{ks}\partial_s\lambda+m \partial_p \lambda=0 \quad \Rightarrow \quad (m+2)\partial_p \lambda=0.
$$
Thus $\lambda$ must be constant. Since $g$ is flat, we can find a change of coordinates  which brings $h$ to constant form. As $L_1$ is a constant scalar operator, it retains its form in any coordinate system.
Similarly $\lambda^i$ and $g_{\lambda^i}$ can be reduced to constant form.
\endproof

This  leads to the following

\begin{theorem}\label{thmMok}
Consider a non-degenerate  Hamiltonian operator  \eqref{Ham} such that the affinor $L^i_j:=\tilde g^{ik} g_{kj}$ has  (pointwise) diagonal  Jordan normal form. Then this  operator can be reduced to  constant coefficient form by a local change of coordinates.
\end{theorem}

This extends the analogous result of Mokhov \cite{Mokhov2} obtained under the additional assumption of simplicity of  the spectrum of $L$.

Suppose that $g$ has Euclidean signature (or, more generally, there exists a non-degenerate Euclidean combination of the form $\lambda g+\mu \tilde g$). Then  the affinor $L$ can be brought to diagonal form. By Theorem \ref{thmMok} we have

\begin{corollary}\label{corMok}
If one of the  contravariant metrics which define a 2D Hamiltonian operator  is Euclidean, then the operator can be reduced to constant coefficient form.
\end{corollary}
This shows that the most interesting case is when each representative of the pencil $\lambda g+\mu \tilde g$ is essentially pseudo-Euclidean, and the affinor  $L$ has non-trivial Jordan block structure.


\section{Classification results}\label{classification}
In this section we classify Hamiltonian operators of type (\ref{Ham}) with the number of components $n\leq 4$. 
This will be done up to arbitrary transformations of the dependent variables $u^i$. 
Our approach  is based on the following two fundamental facts:
\begin{enumerate}
\item Any Killing bivector in flat space is the sum of symmetrized tensor products of Killing vectors;
\item A pair of symmetric bivectors can be brought to the Segre  normal form.
\end{enumerate}
We recall that the first metric $g$ can  always be reduced to constant form, and the second one must be linear, that is
\begin{equation}
 \tilde{g}=c^{ij}_k u^k + g^{ij}_0,
\end{equation}
here $g$ and $g_0$ are constant symmetric matrices, and $c^{ij}_k$ are constant coefficients.
Taking  `generic'  values $u^k_0$ of the variables $u^k$ and applying the shift of variables,
$$
u^k \rightarrow u^k_0+v^k,
$$
we obtain the transformed metric,
\begin{equation}
 \tilde{g}=c^{ij}_k v^k + \tilde{g}^{ij}_{0}.
\end{equation}
The genericity of $u^k_0$  allows us to assume that the Segre type of the pair $(g, \tilde g)$ is the same as that of $(g, \tilde g_0)$. Recall that the Segre type of  a  pair of symmetric forms can be read off the Jordan normal form of the corresponding  affinor $L$, see below. 
Bringing $g$ and $\tilde{g}_{0}$ to the Segre normal form leads to a considerable simplification of  calculations.
Furthermore, the splitting lemma  allows us to consider  irreducible cases only,  where the affinor $L$ either  has  one real eigenvalue, or two complex conjugate eigenvalues. 

The theory of normal forms of pairs of symmetric bilinear forms  is based on the following result, see e.g.   \cite{lancaster}:
\begin{theorem}
Suppose $L$ is a $g$-selfadjoint operator on a real vector space $V$.
There exist a canonical basis $e_1, \ldots, e_n \in V$ in which $L$ and $g$ can be simultaneously reduced to the following block diagonal canonical forms:
$$
L_{can}=
\begin{pmatrix}
L_1 & & & \\
& L_2 & &\\
& & \ddots &\\
& & & L_s
\end{pmatrix},
\quad
g_{can}=
\begin{pmatrix}
g_1 & & & \\
& g_2 & &\\
& & \ddots &\\
& & & g_s
\end{pmatrix},
$$
where
$$
g_{j}=\pm
\begin{pmatrix}
 & & &1 \\
& &1 &\\
&\iddots & &\\
1& & & 
\end{pmatrix},
$$
and
$$
L_{j}=
\begin{pmatrix}
\lambda^j & 1 & & \\
& \lambda^j&\ddots &\\
& &\ddots &1\\
& & & \lambda^j 
\end{pmatrix},
$$
in the case of real eigenvalues $\lambda^j \in \mathbb{R}$ (real Jordan block), or
$$
L_j=
\begin{pmatrix}
 \begin{matrix}
 a & \!\!\! b\\
 -b & a
 \end{matrix}
&
 \begin{matrix}
 1 & 0\\
 0 & 1
 \end{matrix} 
& & \\
&
 \begin{matrix}
 a & \!\!\! b\\
 -b & a
 \end{matrix}
& \ddots &\\
& & \ddots 
 \begin{matrix}
 1 & 0\\
 0 & 1
 \end{matrix}
\\
& & & 
 \begin{matrix}
 a & \!\!\! b\\
 -b & a
 \end{matrix}
\end{pmatrix},
$$
in the case of complex conjugate eigenvalues $\lambda^j_{1,2}=a\pm i b$ (complex Jordan block). It is assumed that for each $j$ the blocks $g_j$ and $L_j$ are of the same size.
\end{theorem}
\begin{remark}
Let us briefly comment on what we mean by Segre type. Suppose $n=4$ and let us consider the affinor $L=\tilde{g} g^{-1}$. 
In the  case of two complex conjugate eigenvalues $\nu+i \lambda$ and  $\nu-i \lambda$, the canonical form of  $L$ reads
$$
\begin{pmatrix}
\nu & -\lambda & 1 & 0\\
\lambda & \nu & 0 & 1 \\
0 & 0 & \nu & -\lambda\\
0 & 0 & \lambda & \nu
\end{pmatrix}.
$$
In the case of a single real eigenvalue we have the following  four canonical forms:
$$
\begin{array}{cccc}
\begin{pmatrix}
\lambda & 1 & 0 & 0\\
0 & \lambda & 1 & 0\\
0 & 0 & \lambda & 1\\
0 & 0 & 0 & \lambda
\end{pmatrix},
&
\begin{pmatrix}
\lambda & 1 & 0 & 0\\
0 & \lambda & 1 & 0\\
0 & 0 & \lambda & 0\\
0 & 0 & 0 & \lambda
\end{pmatrix},
&
\begin{pmatrix}
\lambda & 1 & 0 & 0\\
0 & \lambda & 0 & 0\\
0 & 0 & \lambda & 1\\
0 & 0 & 0 & \lambda
\end{pmatrix},
&
\begin{pmatrix}
\lambda & 1 & 0 & 0\\
0 & \lambda & 0 & 0\\
0 & 0 & \lambda & 0\\
0 & 0 & 0 & \lambda
\end{pmatrix}.
\\[30pt]
\mbox{\footnotesize \it Segre type [4]}  & \mbox{\footnotesize \it Segre type [3,1]} &  \mbox{ \footnotesize\it Segre type [2,2]}  &  \mbox{ \footnotesize\it Segre type [2,1,1]}
\end{array}
$$
Segre type indicates the number and sizes of Jordan blocks with the same eigenvalue $\lambda$.
\end{remark}

\subsection{One-component case}
It was shown in \cite{DN1, Mokhov1} that any one-component operator can be reduced to  constant coefficient form,
$
P=\lambda \partial_x + \mu \partial_y,
$
here $\lambda$ and $\mu$ are arbitrary constants.

\subsection{Two-component case}
The two-component situation is also understood completely \cite{Dubrovin1, Mokhov1}:
we have only one non-constant Hamiltonian operator (\ref{2d_operator}), the corresponding affinor $L$ is a single Jordan block  with non-constant eigenvalue: 
$$
P=
\begin{pmatrix}
0 & 1 \\
1 & 0
\end{pmatrix}
\frac{d}{dx}
+\begin{pmatrix}
-2u^1 & u^2 \\
u^2 & 0
\end{pmatrix}
\frac{d}{dy}
+\begin{pmatrix}
-u_y^1 & 2 u_y^2 \\
-u_y^2 & 0
\end{pmatrix}.
$$
Let us give an alternative proof of this result based on the Killing condition. First we reduce  $g$ to flat coordinates,
$$g=
\begin{pmatrix}
0 & 1 \\
1 & 0
\end{pmatrix},
$$
recall that $g$ must be Lorentzian. Since $\tilde{g}$ is a  Killing tensor of $g$, it is a quadratic expression in the isometries $u^1\partial_{u^1}-u^2 \partial_{u^2}, \partial_{u^1}, \partial_{u^2}$. Since $\tilde g$ is linear, the first isometry can only enter linearly, so that 
$$
\tilde{g}= (u^1 \partial_{u^1} - u^2\partial_{u^2})(\alpha \partial_{u^1} + \beta \partial_{u^2}) + \gamma \partial_{u^1}^2 +2 \delta \partial_{u^1} \partial_{u^2} + \epsilon \partial_{u^2}^2,
$$
here $\alpha, \beta, \gamma, \delta, \epsilon$ are arbitrary constants. The vanishing of the Nijenhuis torsion of the corresponding affinor $L$  gives 
$$
 (\alpha u^1 + \gamma) \beta=0, ~~~
 (\beta u^2 -\epsilon) \alpha=0.
$$
Without any loss of generality one can take $\beta=0$. In this case $\alpha $ must be nonzero, otherwise $\tilde g$ will have constant coefficients. Then $\epsilon=0$, and modulo translations of $u^1, u^2$ we arrive at the required expression \eqref{2d_operator}.

\subsection{Three-component case}
Our main result can be summarised as follows.

\medskip

\noindent{\bf Theorem 3.}
{\it Any  irreducible non-constant three-component Hamiltonian operator in 2D can be brought (by a change of the dependent variables $u^i$) to the form $\pm P$ where $P$ can have one of the following two canonical forms (in both cases the  affinor $L$ is a single $3\times 3$ Jordan block):
\begin{enumerate}
\item
Jordan block with constant eigenvalue
$$
P=
\begin{pmatrix}
0 & 0 & 1 \\
0 & 1 & 0 \\
1 & 0 & 0
\end{pmatrix}
\frac{d}{dx}
+\begin{pmatrix}
-2 u^2 & u^3 & \lambda \\
u^3 & \lambda & 0 \\
\lambda & 0 & 0
\end{pmatrix}
\frac{d}{dy}
+\begin{pmatrix}
-u_y^2 & 2 u_y^3 & 0\\
-u_y^3 & 0 & 0\\
0 & 0 & 0
\end{pmatrix},
$$
\item
Jordan block with non-constant eigenvalue
$$
P=
\begin{pmatrix}
0 & 0 & 1 \\
0 & 1 & 0 \\
1 & 0 & 0
\end{pmatrix}
\frac{d}{dx}
+\begin{pmatrix}
-2u^1 & -\frac{1}{2}u^2 & u^3 \\
-\frac{1}{2}u^2 & u^3 & 0 \\
u^3 & 0 & 0
\end{pmatrix}
\frac{d}{dy}
+\begin{pmatrix}
-u_y^1 &\frac{1}{2} u_y^2 & 2u_y^3\\
-u_y^2 & \frac{1}{2} u_y^3 & 0\\
-u_y^3 & 0 & 0
\end{pmatrix}.
$$
\end{enumerate}
}

\proof

Since the complex conjugate case cannot occur (it requires an even number of components), we only need to consider the cases where the affinor $L$ has one triple eigenvalue, and has Segre type $[3]$ or $[2, 1]$. Since the case  $[2, 1]$ gives no non-constant examples, we will concentrate on  Segre type $[3]$. Then there exists a coordinate system where  $g$ and $\tilde{g}_0$ take the form
$$
g^{ij}=
\begin{pmatrix}
 0 & 0 & 1 \\
0 & 1 & 0 \\
 1 & 0 & 0 
\end{pmatrix},
\quad
\tilde{g}^{ij}_{0}=
\begin{pmatrix}
 0 & 1 & \lambda \\
 1 & \lambda & 0 \\
 \lambda & 0 & 0 
\end{pmatrix}.
$$
The general solution of Mokhov's conditions is given by the two-parameter family $\tilde{g}= \kappa_1 \tilde{g}_1+\kappa_2 \tilde{g}_2+\tilde{g}_{0}$, where $\kappa_i$ are arbitrary constants, and the  bivectors $\tilde g_i$ are as follows:
$$
\tilde{g}_1=
\begin{pmatrix}
-2 u^1 & -\frac{1}{2} u^2 & u^3 \\
-\frac{1}{2} u^2 & u^3 & 0 \\
u^3 & 0 & 0 
\end{pmatrix},
\quad
\tilde{g}_2=
\begin{pmatrix}
- 2u^2 & u^3 & 0 \\
u^3 & 0 & 0 \\
0 & 0 & 0
\end{pmatrix}.
$$
In the non-constant eigenvalue case,  $\kappa_1\neq0$, using the following transformations which preserve both $g$ and $\tilde{g}_0$, 
$$
\begin{array}{ccl}
u^1 &\to& \dfrac{1}{\kappa_1} u^1+\dfrac{2\kappa_2}{\kappa_1^2} u^2-\dfrac{2\kappa_2^2}{\kappa_1^3} u^3-\dfrac{2\kappa_2}{\kappa_1^3},\\[5pt]
u^2 &\to& u^2-\dfrac{2\kappa_2}{\kappa_1} u^3 +\dfrac{2 \kappa_1-2}{\kappa_1},\\[5pt]
u^3 &\to& \kappa_1 u^3,
\end{array}
$$
we can reduce the above family  to  $\tilde{g}=\tilde{g}_1+\tilde{g}_0$ (that is, we can set $\kappa_1=1, \ \kappa_2=0$).  After that we can eliminate $\tilde g_0$ by  appropriate translations of   $u^2$ and $u^3$, arriving at the final answer $\tilde{g}=\tilde{g}_1$.
Similarly, in the constant eigenvalue case,  $\kappa_1=0$, we can set  $\kappa_2=1$, and use an  appropriate translations  of $u^3$ to arrive at the normal form above.

 
\endproof

\subsection{Four-component case}\label{sect_4}
The four-component situation is more complicated since we have more Segre types. In this section we present the results of classification  of four-component Hamiltonian operators of the form \eqref{Ham} with one real eigenvalue, as well as with  two complex conjugate eigenvalues (the latter turn out to be complexifications of  the $2\times 2$ operator \eqref{2d_operator}).
We will only give  canonical forms for the contravariant metrics $g, \tilde g$: the symbols $\tilde{b}^{ij}_k$ of the second metric can be computed directly. We skip the details of calculations: these follow the procedure outlined at the beginning of Section 5, and are essentially the same as in the proof of Theorem 3. 


\subsubsection{Segre type [2,1,1]}

One can show that this case leads to constant coefficient operators.

\subsubsection{Segre type [2,2]}
By Theorem 7, we have to consider two different cases.

\noindent
{\bf Case 1}:
There exists a coordinate system where  $g$ and $\tilde{g}_0$ take the form
$$
g^{ij}=
\begin{pmatrix}
0 & 1 & 0 & 0\\
1 & 0 & 0 & 0 \\
0 & 0 & 0 & 1\\
0 & 0 & 1 & 0
\end{pmatrix},
\quad
\tilde{g}^{ij}_{0}=
\begin{pmatrix}
1 & \lambda & 0 & 0\\
\lambda & 0 & 0 & 0 \\
0 & 0 & 1 &  \lambda \\
0 & 0 &  \lambda & 0
\end{pmatrix}.
$$
The general solution of Mokhov's conditions is given by $\tilde{g}=\sum_{i=1}^{4} \kappa_i \tilde{g}_i+\tilde{g}_{0}$, where $\kappa_i$ are arbitrary constants, and the  bivectors $\tilde g_i$ are as follows:
$$
\tilde{g}_1=
\begin{pmatrix}
u^1 & -\frac{1}{2} u^2 & \frac{1}{2} u^3 & 0\\
-\frac{1}{2} u^2 & 0 & 0 & 0\\
\frac{1}{2} u^3 & 0 & 0 & - \frac{1}{2} u^2\\
0 & 0 & - \frac{1}{2} u^2 & 0
\end{pmatrix}
,\quad
\tilde{g}_2=
\begin{pmatrix}
u^4 & 0 & -\frac{1}{2} u^2 & 0\\
0 & 0 & 0 & 0\\
-\frac{1}{2} u^2& 0 & 0 & 0\\
0 & 0 & 0 & 0
\end{pmatrix},
$$
$$
\tilde{g}_3=
\begin{pmatrix}
0 & \frac{1}{2} u^4 & - \frac{1}{2} u^1 & 0\\
\frac{1}{2} u^4 & 0 & 0 & 0\\
- \frac{1}{2} u^1 & 0 & -u^3 &  \frac{1}{2} u^4\\
0 & 0 & \frac{1}{2} u^4 & 0
\end{pmatrix}
,\quad
\tilde{g}_4=
\begin{pmatrix}
0& 0 & \frac{1}{2} u^4 & 0\\
0 & 0 & 0 & 0\\
\frac{1}{2} u^4& 0 & - u^2 & 0\\
0 & 0 & 0 & 0
\end{pmatrix}.
$$
The eigenvalue of the corresponding affinor $L$ is $\frac{1}{2}(\kappa_3 u^4-\kappa_1 u^2) +\lambda$.
 Using symmetries which preserve $g$ and $\tilde{g}_{0}$ one can set the coefficients $\kappa_3$ and $\kappa_4$ equal to zero, arriving at  the normal form 
$$
\tilde{g}=\kappa_1 \tilde{g}_1+\kappa_2 \tilde{g}_2+\tilde{g}_{0}.
$$

\noindent
{\bf Case 2}: There exists a coordinate system where  $g$ and $\tilde{g}_0$ take the form
$$
g^{ij}=
\begin{pmatrix}
0 & 1 & 0 & 0\\
1 & 0 & 0 & 0 \\
0 & 0 & 0 & -1\\
0 & 0 & -1 & 0
\end{pmatrix},
\quad
\tilde{g}^{ij}_{0}=
\begin{pmatrix}
1 & \lambda & 0 & 0\\
\lambda & 0 & 0 & 0 \\
0 & 0 & -1 & - \lambda \\
0 & 0 & - \lambda & 0
\end{pmatrix}.
$$
The eigenvalue of the corresponding affinor $L$ is $\frac{1}{2}(\kappa_3 u^4-\kappa_1 u^2) +\lambda$. The general solution of Mokhov's conditions is given by $\tilde{g}=\sum_{i=1}^{4} \kappa_i \tilde{g}_i+\tilde{g}_{0}$, where $\kappa_i$ are arbitrary constants,  and the  bivectors $\tilde g_i$ are as follows:
$$
\tilde{g}_1=
\begin{pmatrix}
u^1 & -\frac{1}{2} u^2 & \frac{1}{2} u^3 & 0\\
-\frac{1}{2} u^2 & 0 & 0 & 0\\
\frac{1}{2} u^3 & 0 & 0 &  \frac{1}{2} u^2\\
0 & 0 &  \frac{1}{2} u^2 & 0
\end{pmatrix}
,\quad
\tilde{g}_2=
\begin{pmatrix}
u^4 & 0 & \frac{1}{2} u^2 & 0\\
0 & 0 & 0 & 0\\
\frac{1}{2} u^2& 0 & 0 & 0\\
0 & 0 & 0 & 0
\end{pmatrix},
$$
$$
\tilde{g}_3=
\begin{pmatrix}
0 & \frac{1}{2} u^4 &  \frac{1}{2} u^1 & 0\\
\frac{1}{2} u^4 & 0 & 0 & 0\\
 \frac{1}{2} u^1 & 0 &  u^3 &- \frac{1}{2} u^4\\
0 & 0 & -\frac{1}{2} u^4 & 0
\end{pmatrix}
,\quad
\tilde{g}_4=
\begin{pmatrix}
0& 0 & \frac{1}{2} u^4 & 0\\
0 & 0 & 0 & 0\\
\frac{1}{2} u^4& 0 &  u^2 & 0\\
0 & 0 & 0 & 0
\end{pmatrix}.
$$
Using symmetries which preserve $g$ and $\tilde{g}_{0}$ one can reduce the above four-parameter family to one of the following  normal forms:
$$
\tilde{g}=\tilde{g}_{0}+
\left\{
\begin{array}{l}
\kappa_1 \tilde{g}_1+\kappa_2 \tilde{g}_4\\
\kappa_1 \tilde{g}_2+\kappa_2 \tilde{g}_3\\
\tilde{g}_2\pm \tilde{g}_4\\
\tilde{g}_1\pm \tilde{g}_3+\kappa_1 \tilde{g}_4
\end{array}
\right.
\quad \quad
\kappa_1,\kappa_2=const.
$$

\subsubsection{Segre type [3,1]}
Here we also have  two different cases.

\noindent
{\bf Case 1}:
There exists a coordinate system where $g$ and $\tilde{g}_{0}$ take the form
$$
g^{ij}=\begin{pmatrix}
0 & 0 & 1 & 0\\
0 & 1 & 0 & 0\\
1 & 0 & 0 & 0\\
0 & 0 & 0 & 1
\end{pmatrix},
\quad
\tilde{g}^{ij}_{0}=\begin{pmatrix}
0 & 1 & \lambda & 0\\
1 & \lambda & 0 & 0\\
\lambda & 0 & 0 & 0\\
0 & 0 & 0 & \lambda
\end{pmatrix}.
$$
The general solution of  Mokhov's conditions is given by $\tilde{g}=\sum_{i=1}^{4} \kappa_i \tilde{g}_i+\tilde{g}_{0}$, where $\kappa_i$ are arbitrary constants, and the  bivectors $\tilde g_i$ are as follows:
$$
\tilde{g}_1=
\begin{pmatrix}
2u^1 & \frac{1}{2} u^2 &-u^3 &  \frac{1}{2} u^4\\
\frac{1}{2} u^2 & -u^3 & 0 & 0\\
-u^3 & 0 & 0 & 0\\
 \frac{1}{2} u^4 & 0 & 0 & - u^3
\end{pmatrix}
,\quad
\tilde{g}_2=
\begin{pmatrix}
u^2 & -\frac{1}{2} u^3 & 0 & 0\\
-\frac{1}{2} u^3 & 0 & 0 & 0\\
0 & 0 & 0 & 0\\
0 & 0 & 0& 0
\end{pmatrix},
$$
$$
\tilde{g}_3=
\begin{pmatrix}
u^4 & 0 & 0 & - \frac{1}{2} u^3\\
0 & 0 & 0 & 0\\
0 & 0 & 0 & 0\\
- \frac{1}{2} u^3 & 0 & 0& 0
\end{pmatrix}
,\quad
\tilde{g}_4=
\begin{pmatrix}
0& \frac{1}{2} u^4 & 0 & - \frac{1}{2} u^2\\
\frac{1}{2} u^4 & 0 & 0 & 0\\
0 & 0 & 0 & 0\\
- \frac{1}{2} u^2 & 0 & 0 & 0
\end{pmatrix}.
$$
The eigenvalue of the corresponding affinor $L$ is $\lambda-\kappa_1 u^3$.
Using symmetries which preserve  $g$ and $\tilde{g}_0$ one can bring the above four-parameter family  to one of the following canonical forms:
$$
\tilde{g}=
\tilde{g}_{0}+
\left\{
\begin{array}{l}
\kappa_1 \tilde{g}_2+\kappa_2 \tilde{g}_3\\
\kappa_1 \tilde{g}_3+\kappa_2 \tilde{g}_4\\
\kappa_1 \tilde{g}_1+\kappa_2 \tilde{g}_2+\kappa_3 \tilde{g}_4
\end{array}
\right.
\quad \quad
\kappa_1,\kappa_2,\kappa_3=const.
$$

\noindent
{\bf Case 2}:
There exists a coordinate system where $g$ and $\tilde{g}_{0}$ take the form
$$
g^{ij}=\begin{pmatrix}
0 & 0 & 1 & 0\\
0 & 1 & 0 & 0\\
1 & 0 & 0 & 0\\
0 & 0 & 0 & -1
\end{pmatrix},
\quad
\tilde{g}^{ij}_{0}=\begin{pmatrix}
0 & 1 & \lambda & 0\\
1 & \lambda & 0 & 0\\
\lambda & 0 & 0 & 0\\
0 & 0 & 0 &- \lambda
\end{pmatrix}.
$$
The general solution of Mokhov's conditions is given by $\tilde{g}=\sum_{i=1}^{4} \kappa_i \tilde{g}_i+\tilde{g}_{0}$, where $\kappa_i$ are arbitrary constants, and the  bivectors $\tilde g_i$ are as follows:
$$
\tilde{g}_1=
\begin{pmatrix}
2u^1 & \frac{1}{2} u^2 &-u^3 &  \frac{1}{2} u^4\\
\frac{1}{2} u^2 & -u^3 & 0 & 0\\
-u^3 & 0 & 0 & 0\\
 \frac{1}{2} u^4 & 0 & 0 &  u^3
\end{pmatrix}
,\quad
\tilde{g}_2=
\begin{pmatrix}
u^2 & -\frac{1}{2} u^3 & 0 & 0\\
-\frac{1}{2} u^3 & 0 & 0 & 0\\
0 & 0 & 0 & 0\\
0 & 0 & 0& 0
\end{pmatrix},
$$
$$
\tilde{g}_3=
\begin{pmatrix}
u^4 & 0 & 0 &  \frac{1}{2} u^3\\
0 & 0 & 0 & 0\\
0 & 0 & 0 & 0\\
 \frac{1}{2} u^3 & 0 & 0& 0
\end{pmatrix}
,\quad
\tilde{g}_4=
\begin{pmatrix}
0& \frac{1}{2} u^4 & 0 &  \frac{1}{2} u^2\\
\frac{1}{2} u^4 & 0 & 0 & 0\\
0 & 0 & 0 & 0\\
 \frac{1}{2} u^2 & 0 & 0 & 0
\end{pmatrix}.
$$
The eigenvalue of the corresponding affinor $L$ is $\lambda-\kappa_1 u^3$.
Using symmetries which preserve $g$ and $\tilde{g}_0$ one can bring the above four-parameter family to one of the following normal forms:
$$
\tilde{g}=
\tilde{g}_{0}+
\left\{
\begin{array}{l}
\kappa_1 \tilde{g}_2+\kappa_2 \tilde{g}_3\\
\kappa_1 \tilde{g}_3+\kappa_2 \tilde{g}_4\\
\kappa_1 \tilde{g}_1+\kappa_2 \tilde{g}_2+\kappa_3 \tilde{g}_4
\end{array}
\right.
\quad \quad
\kappa_1,\kappa_2,\kappa_3=const.
$$

\subsubsection{Segre type [4]}
This is the case where the corresponding affinor $L$ is a single Jordan block (see Section 6 for the general theory). There exists a coordinate system where  $g$ and $\tilde{g}_{0}$ take the form
$$
g^{ij}=
\begin{pmatrix}
0 & 0 & 0 & 1 \\
0 & 0 & 1 & 0 \\
0 & 1 & 0 & 0 \\
1 & 0 & 0 & 0
\end{pmatrix},
\quad
\tilde{g}^{ij}_{0}=
\begin{pmatrix}
0 & 0 & 1 & \lambda \\
0 & 1 & \lambda & 0 \\
1 & \lambda & 0 & 0 \\
\lambda & 0 & 0 & 0
\end{pmatrix}.
$$
It turns out that the general solution of Mokhov's conditions  is
$\tilde{g}=\sum_{i=1}^{3} \kappa_i \tilde{g}^i+\tilde{g}_{0}$
where $\kappa_i$ are arbitrary constants, and the  bivectors $\tilde g_i$ are as follows:
$$
\tilde{g}_1=
\begin{pmatrix}
-u^1 & -\frac{1}{2} u^2 & 0 & \frac{1}{2} u^4\\
-\frac{1}{2} u^2 & 0 & \frac{1}{2} u^4& 0 \\
0 & \frac{1}{2} u^4& 0 &0 \\
\frac{1}{2} u^4& 0 &0 &0
\end{pmatrix},
\tilde{g}_2=
\begin{pmatrix}
2 u^2 & \frac{1}{2} u^3 & -u^4 &0\\
\frac{1}{2} u^3 & -u^4 & 0& 0 \\
-u^4 & 0 & 0 &0 \\
0& 0 &0 &0
\end{pmatrix},
\tilde{g}_3=
\begin{pmatrix}
u^3 & -\frac{1}{2} u^4 & 0 & 0\\
-\frac{1}{2} u^4 & 0 & 0& 0 \\
0 & 0& 0 &0 \\
0& 0 &0 &0
\end{pmatrix}.
$$
Here the eigenvalue of the affinor $L$ is $\frac{1}{2}\kappa_1 u^4+\lambda$.
Using symmetries which preserve $g$ and $\tilde{g}_0$ one can bring  the above three-parameter family to one of the following normal forms:
in the non-constant eigenvalue case
$$
\tilde{g}=\tilde{g}_0+\tilde{g}_1+\kappa_1\tilde{g}_2, \quad \kappa_1=const,
$$
while in the constant eigenvalue case
$$
\tilde{g}=\tilde{g}_0+
\left\{
\begin{array}{l}
\tilde{g}_2\\
k_1 \tilde{g}_3
\end{array}
\right.
\quad \kappa_1=const.
$$

\subsubsection{Complex conjugate case}
In the case of two pairs of complex conjugate eigenvalues $\nu + i \lambda$ and $\nu - i \lambda$, there exists a coordinate system such that
$$
g^{ij}=
\begin{pmatrix}
0 & 0 & 0 & 1\\
0 & 0 & 1 & 0 \\
0 & 1 & 0 & 0\\
1 & 0 & 0 & 0
\end{pmatrix},
\quad
\tilde{g}^{ij}_{0}=
\begin{pmatrix}
0 & 1 & -\lambda & \nu\\
1 & 0 & \nu & \lambda \\
-\lambda & \nu & 0 & 0\\
\nu & \lambda & 0 & 0
\end{pmatrix}.
$$
The general solution of Mokhov's conditions is $\tilde{g}=\kappa_1\tilde{g}_1+\kappa_2\tilde{g}_2+\tilde{g}_0$ where  $\kappa_i$ are arbitrary constants, and  the  bivectors $\tilde g_i$ are as follows:
$$
\tilde{g}_1=
\begin{pmatrix}
2 u^2 & - 2u^1 &  -u^4 & u^3\\
- 2u^1& -  2u^2 &  u^3 & u^4 \\
 - u^4 & u^3 & 0 & 0\\
u^3 &  u^4 & 0 & 0
\end{pmatrix},
\quad
\tilde{g}_2=
\begin{pmatrix}
2 u^1 & 2 u^2 & -u^3 & -u^4 \\
2 u^2 & - 2 u^1&  - u^4 &  u^3 \\
- u^3 & - u^4 & 0 & 0\\
-u^4 &  u^3  & 0 & 0
\end{pmatrix}.
$$
Using symmetries which preserve the form of $g$  one can eliminate $\tilde g_0$, and  bring $\tilde g$ to the normal form
$$
\tilde{g}^{ij}=
\begin{pmatrix}
2 u^2 & -2 u^1 & - u^4 & u^3\\
-2 u^1 & -2 u^2 & u^3 & u^4 \\
-  u^4 & u^3 & 0 & 0\\
u^3 & u^4 & 0 & 0
\end{pmatrix}.
$$
The eigenvalues of the corresponding affinor $L$ are $u^3\pm i u^4$.
Note  that this case is a complexification of the two-component operator (\ref{2d_operator}), which can be achieved via the following recipe (see \cite{BM} for more details):
each complex entry $a+ib$ of $g^{\mathbb{C}}$ and $\tilde{g}^{\mathbb{C}}$ is replaced by the $2 \times 2$ block
$$
\begin{pmatrix}
-b & a\\
a & b
\end{pmatrix},
$$
where $g^{\mathbb{C}}$ and $\tilde{g}^{\mathbb{C}}$ are the complexified bivectors of the operator  (\ref{2d_operator}):
$$
g^{\mathbb{C}}=
\begin{pmatrix}
0 & 1\\
1 & 0
\end{pmatrix},
\quad
\tilde{g}^{\mathbb{C}}=
\begin{pmatrix}
-2 z^1 & z^2\\
z^2 & 0
\end{pmatrix},
\quad
$$
 $z^1=u^1+iu^2,\ z^2=u^3+iu^4.$


\section{The single Jordan block case}\label{sect_Jordan}

Let us begin with  examples of $n$-component Hamiltonian operators of the single Jordan block type.

\medskip

\noindent {\bf Example 1.}
One of the most important examples  was discovered by Mokhov \cite{Mokhov1}. Here the first $n \times n$ contravariant metric is constant and anti-diagonal,
$$
g=\begin{pmatrix}
& & 1\\
& \iddots &\\
1 & & 
\end{pmatrix},
$$
while the second contravariant metric $\tilde g$ is defined as follows:
$$
\tilde g^{ij}=(b^{ij}_k+b^{ji}_k) u^k, \quad
\left\{\begin{array}{ll}
b^{ij}_k=0 & \mbox{if } k\neq i+j-1,\\
b^{ij}_{i+j-1}=3j-n-2 & \mbox{otherwise.}
\end{array}
\right.
$$
One can verify that  the Jordan normal form for the corresponding affinor $L$ is a single Jordan block with non-constant eigenvalue (for any $n \neq 4$: in the exceptional case $n=4$ the affinor $L$ is the sum of two $2\times 2$ Jordan blocks). We will  refer to this case as the Mokhov operator. The affinor $L$ is given by  $L^i_j=[3(i-j)+n-1]u^{n+i-j}$.  The equivalent form for $\tilde g$ is
 $$
 \tilde g^{ij}=[3(i+j)-2(n+2)]u^{i+j-1},
 $$
 for $i+j-1 \le n$, and $0$ otherwise (in what follows, we use the following convention: if $\alpha>n$ then $u^{\alpha} \equiv 0$).
For  $n=2,3,4$ the explicit form of $\tilde g$ is as follows:
$$
\tilde g=
\begin{pmatrix}
-2u^1& u^2\\
u^2 & 0
\end{pmatrix},
\quad
\tilde g=
\begin{pmatrix}
-4u^1& -u^2 & 2 u^3\\
-u^2 & 2 u^3 & 0\\
2u^3 & 0 & 0
\end{pmatrix},
\quad
\tilde g=
\begin{pmatrix}
-6u^1& -3u^2 &0 & 3 u^3\\
-3u^2 &0& 3 u^3 & 0\\
0 & 3 u^3 & 0 & 0\\
3 u^3 & 0 & 0 & 0
\end{pmatrix}.
$$

\medskip

\noindent {\bf Example 2.}
Another  $n$-component example has $g$ the same as in Example 1, while the second contravariant metric is
$$
\tilde{g}^{ij}=(b^{ij}_k+b^{ji}_k) u^k+\lambda g^{ij}, \quad
\left\{\begin{array}{ll}
b^{ij}_k=0 & \mbox{if } k\neq i+j,\\
b^{ij}_{i+j}=3j-n-1 & \mbox{otherwise,}
\end{array}
\right.
$$
$\lambda=const$. One can verify that this pair of contravariant metrics defines a Hamiltonian operator for any $ n \ge3$ (the case $n=2$ is trivial since all  $b^{ij}_k$ vanish). The corresponding affinor $L$ is a single Jordan block with constant eigenvalue $\lambda$.
For instance, for $n=3,4$ the second contravariant metric reads 
$$
\tilde{g}=
\begin{pmatrix}
-2u^2& u^3 &  \lambda\\
u^3 &  \lambda & 0\\
\lambda & 0 & 0
\end{pmatrix},
\quad
\tilde{g}=
\begin{pmatrix}
-4u^2& -u^3 & 2 u^4 &  \lambda\\
-u^3 & 2 u^4 & \lambda& 0\\
2 u^4 &  \lambda & 0 & 0\\
\lambda & 0 & 0 & 0
\end{pmatrix}.
$$

\medskip

The aim of this section is to give a complete description of the case where the affinor $L$ is a single Jordan block. We will see that the Mokhov example plays fundamental role in this picture. 
To formulate our main result, let us introduce  symmetric bivectors  $\mu^{(n;k)}$ as follows:
\begin{equation}\label{Mok_shift}
\mu^{(n;k) ij}=[3(i+j)-2(n+2-k)]u^{i+j-1+k}.
\end{equation}
In particular, $\mu^{(n;0)}$ coincides with the second contravariant metric $\tilde g$ of the   Mokhov operator from Example 1. Note also that $\mu^{(n;k)}=0$ for $k>n-2$.
Let us present the explicit form for some  $\mu^{(n;k)}$:
$$
\mu^{(3;1)}=
\begin{pmatrix}
-2u^2& u^3 &  0\\
u^3 &  0 & 0\\
0 & 0 & 0
\end{pmatrix},
\quad
\mu^{(4;1)}=
\begin{pmatrix}
-4u^2& -u^3 & 2 u^4 &  0\\
-u^3 & 2 u^4 & 0& 0\\
2 u^4 &  0 & 0 & 0\\
0 & 0 & 0 & 0
\end{pmatrix},
$$
$$
\mu^{(4;2)}=
\begin{pmatrix}
-2u^3& u^4 & 0 &  0\\
u^4 & 0 & 0& 0\\
0 &  0 & 0 & 0\\
0 & 0 & 0 & 0
\end{pmatrix}.
$$

We will show that in the case when the affinor $L$ is a single Jordan block,  the general solution of Mokhov's conditions reads
\begin{equation}
 \tilde{g}=\tilde{g}_{0}+\sum_{m=0}^{n-2} \xi_m \mu^{(n;m)},
\label{Jor_sol}
\end{equation}
where $\xi_m$ are arbitrary constants, and
$$
g= \pm
\begin{pmatrix}
& &&  1\\
& & 1 &\\
& \iddots & & \\
1 & & &
\end{pmatrix},
\quad
\tilde{g}_{0}= \pm
\begin{pmatrix}
& &1 & \lambda\\
& \iddots & \lambda& \\
1& \iddots&  & \\
\lambda & & &
\end{pmatrix}.
$$
Here the eigenvalue of $L$ equals $\xi_0 (n-1)u^n+\lambda$. In the non-constant eigenvalue  case, $\xi_0 \ne 0$, we have the following result:

\medskip

\vspace{5pt}\noindent
{\bf Theorem \ref{Jor_thm}} 
{\it Let $P$ be a Hamiltonian operator (\ref{Ham}) such that the affinor  $L= \tilde{g} g^{-1}$ is a single $n\times n$ Jordan block with non-constant eigenvalue. Then there exists a coordinate system in which $g$  and $\tilde g$ can be reduced to the following canonical forms:
$$
g=\pm
\begin{pmatrix}
 &&  1\\
& \iddots &  \\
1 & & 
\end{pmatrix},
\quad
\tilde{g}=\pm
\left\{
\begin{array}{lcl}
\mu^{(n;0)} & \mbox{if} & n\not\equiv 1\!\!\!\! \mod 3,\\
\mu^{(n;0)} + \kappa_1 \mu^{(n;\frac{n-1}{3})} & \mbox{if} & n\equiv 1\!\!\!\! \mod 3, n \neq 4,\\
\mu^{(4;0)}+\kappa_1 \mu^{(4;1)}+ \tilde{g}_0 & \mbox{if} & n=4.
\end{array}
\right.
$$
Here $ \kappa_1$ is an arbitrary constant.  
}

\vspace{5pt}

In the constant eigenvalue case, $\xi_0=0$, we have several canonical forms depending on how many coefficients among $\xi_i$ are equal to zero:
\begin{theorem}[Constant eigenvalue case]\label{Jor_thm_2}
Suppose $\xi_i=0$ for $i=0, \ldots, \alpha-1$. Then the family (\ref{Jor_sol})  can be reduced to
$$
\tilde{g}=\mu^{(n;\alpha)}+\kappa \mu^{(n;\alpha+m)}+\tilde{g}_0, \quad m=\frac{n-1+2\alpha}{3},
$$
if $m \in \mathbb{N}$, otherwise to
$$
\tilde g=\mu^{(n;\alpha)}+\tilde{g}_0.
$$
\end{theorem}

\subsection{Proof of Theorem \ref{Jor_thm}}

The idea of the proof is as follows: first, we find the general solution of Mokhov's equations. It turns out (Proposition 3) that this solution depends on $n-1$ parameters. Using orthogonal transformations, we then reduce this $(n-1)$-parameter family  to various normal forms (Lemma \ref{Norm1} and Proposition  \ref{Norm2}).
We will work in coordinates where  $g$ and $\tilde{g}_{0}$ take canonical  form
$$
g= \pm
\begin{pmatrix}
& &&  1\\
& & 1 &\\
& \iddots & & \\
1 & & &
\end{pmatrix},
\quad
\tilde{g}_{0}= \pm
\begin{pmatrix}
& &1 & \lambda\\
& \iddots & \lambda& \\
1& \iddots&  & \\
\lambda & & &
\end{pmatrix}.
$$
For definiteness, we will consider the $+$ sign. In what follows we will need the following result: 
\begin{proposition}\label{IsoN}
The Killing vectors of $g$ are the following $\frac{1}{2}n(n-1)$ vector fields:
$$
X_{(\alpha,\beta)}=u^{\alpha} \partial_{\beta}-u^{n+1-\beta} \partial_{n+1-\alpha}, \quad X_{\gamma}=\partial_{\gamma},
$$
here $\alpha +\beta < n+1$, and $\partial_{\alpha}=\partial_{u^{\alpha}}$.
\end{proposition}


The  affinor $L$  and the metric $\tilde{g}$ are given by
\begin{eqnarray*}
L^{i}_{j}&=&c^i_{jk}u^k+\tilde{g}^{il}_{0}g_{lj},\\
\tilde{g}^{ij}&=&L^i_lg^{lj}=c^i_{n+1-j,k}u^k+\tilde{g}^{ij}_{0}.
\end{eqnarray*}
These have to satisfy a set of constraints (note that the vanishing of the Nijenhuis torsion, $\mathcal{N}(L)=0$, gives  two types of relations: linear and quadratic in $c^i_{jk}$):
\begin{itemize}
\item Linear part of the condition $\mathcal{N}(L)=0$ reads
\begin{eqnarray}
\label{linN}
c^k_{j,i-1}-c^k_{i,j-1}+c^{k+1}_{ij}-c^{k+1}_{ji}=0;
\end{eqnarray}
\item Quadratic part of the condition $\mathcal{N}(L)=0$ reads
$$c^s_{il}c^m_{js}-c^s_{jl}c^m_{is}+c^m_{sl}c^s_{ij}-c^m_{sl}c^s_{ji}=0;$$
\item Symmetry of $\tilde{g}$ gives
\begin{equation}\label{symmetry}
c^{n+1-i}_{jk}=c^{n+1-j}_{ik};
\end{equation}
\item The Killing condition gives
\begin{equation}\label{killing_jordan}
c^{n+1-i}_{jk}+c^{n+1-k}_{ij}+c^{n+1-j}_{ki}=0.
\end{equation}
\end{itemize}

Remarkably, the linear system (\ref{linN})-(\ref{killing_jordan}) can be solved explicitly:

\begin{proposition}
The general solution of the linear system  (\ref{linN})-(\ref{killing_jordan}) is given by (\ref{Jor_sol}), 
$$
\tilde{g}=\tilde{g}_{0}+\sum_{m=0}^{n-2} \xi_m \mu^{(n;m)},
$$
where $\xi_m$ are arbitrary constants. The eigenvalue of the corresponding  affinor $L$ is $\xi_0 (n-1)u^n+\lambda $.
\end{proposition}

\proof The key  observation allowing one to prove Proposition 3 by induction is as follows.
Suppose  $c^n_{ji}=0$. In this case it is easy to see that  $c^k_{1i}$ and $c^k_{j1}$ must also vanish,
indeed, from \eqref{symmetry} we have
\begin{equation*}
c^{n+1-j}_{1k}=c^{n}_{jk},
\end{equation*}    
and from \eqref{killing_jordan} and \eqref{symmetry} we obtain
\begin{equation*}
c^{n+1-j}_{1k}+c^{n+1-k}_{1j}+c^{n+1-j}_{k1}=0.
\end{equation*}
Then the remaining equations for $c^{k}_{ij}$, with $i,j=2,...,n$ and $k=1,...,n-1$, coincide with the system one obtains in the $(n-1)$-component case with $\tilde{c}^{k}_{ij}=c^{k}_{i-1,j-1}$,  allowing one to use inductive assumption. 



Our strategy will be the following: first we show that the above equations imply  $c^n_{ji}=0$, $c^k_{1i}=0$ and $c^k_{j1}=0$ apart from $c^n_{nn}=c$, $c^1_{1n}=c$ and $c^1_{n1}=-2c$. We already know a solution with $c\ne 0$, which comes from  Mokhov's example. The generic  solution can be written as a linear combination of Mokhov's solution and a solution of the system with $c^n_{ji}=c^k_{1i}=c^k_{j1}=0$ (in which case we can use inductive assumption as outlined above).
 
 From \eqref{symmetry} and \eqref{linN} we have (for $j\ne 1$)
 $$c^k_{1,j-1}=c^n_{n+1-k,j-1}=c^n_{j,n-k}.$$
 Using this identity we can write \eqref{linN} as
 $$c^n_{j,n-k}-c^{k+1}_{1j}+c^{k+1}_{j1}=0.$$
 Similarly,  from \eqref{killing_jordan} we obtain
 $$c^n_{j,n-k}+c^{k+1}_{j1}+c^{k+1}_{1j}=0.$$
 Combining these two conditions we get $c^i_{1j}=0,$ for any $ i\ne 1,j\ne 1$. 
Writing out \eqref{killing_jordan} with $j=k=1$ we get
 $$c^{n+1-i}_{11}+c^{n}_{1i}+c^{n}_{i1}=2c^{n+1-i}_{11}+c^{n}_{1i}=0,$$
which implies $c^{k}_{11}=0$.
  Summarizing, we have
 $$c^i_{1j}=0,\qquad\forall i\ne 1,$$
 which, for symmetry reasons, implies
 $$c^{n}_{jk}=0,\qquad\forall j\ne n,$$
 and
 $$c^i_{j1}=0,\qquad\forall i\ne 1,j\ne n.$$
Our next remark is that $c^1_{1j}=0$ for $j=1,...,n-2$. This  follows from \eqref{linN} evaluated at $k=i=1,$
$$c^1_{1,j-1}-c^{2}_{1j}+c^{2}_{j1}=0.$$
This  readily implies  $c^n_{n,k}=0$ and $c^1_{k,1}=0$ for $k=1,\dots,n-2$, as well as $c^k_{n,1}=0$ for $k=3,\dots,n$.
It is also easy to see that the three non-vanishing coefficients
 $c^n_{nn},c^1_{1n},c^1_{n1}$ are related by 
$$c^n_{nn}=c^1_{1n},\qquad c^1_{n1}=-2c^1_{1n}.$$
We  still need to prove that $c^1_{1,n-1}=0$. Due to the above computations 
 the first column of  the affinor $L$  has the form $(\nu, 0, \dots, 0)^t$ where $\nu=c^1_{1,n-1}u^{n-1}+c^1_{1,n}u^{n}+\lambda$ is the (unique) eigenvalue of $L$. 
Similarly, the last row of $L$ is given by $(0,\dots,0,\nu)$. 
Let us denote by $(e_{(1)},\dots,e_{(n)})$ the canonical frame of the pair $(L, g)$. Thus, 
\begin{eqnarray}
\label{cf1}
L^i_k e^k_{(p)}&=&\nu e^i_{(p)}+e^i_{(p-1)}.
\end{eqnarray}
It follows from the vanishing of the Nijenhuis torsion of $L$ that $e_{(i)}(\nu)=0$ for $i=1,\ldots, n-1$, where $e_{(i)}(\nu)$ denotes the Lie derivative of $\nu$ in direction $e_{(i)}$, see \cite{BM}. Due to the form of the affinor we have (set $i=n$ in \eqref{cf1}):
$$e^n_{(p)}=0,\qquad p=1,\dots,n-1.$$
This means that $e_{(1)},\dots,e_{(n-1)}$ do not contain $\frac{\d}{\d u^n}$,
 and thus $\nu$ must depend  on $u^n$ only, so that   $c^1_{1,n-1}=0$. 

This proves that the general solution  is given by \eqref{Jor_sol}. A direct computations shows  that \eqref{Jor_sol} also satisfies  the quadratic conditions coming from $\mathcal{N}(L)=0$.
\endproof

Thus, the general solution depends on $n-1$ parameters. At this point one might wonder whether this number can be reduced. The answer is yes,  the list of normal forms is presented below.
In order to proceed, we need the following statement.

\begin{lemma}\label{lemmaLie}
The $n-2$ vector fields
$
X_{(k)}=\sum_{i=1}^{n-k}(n-k+1-2i) u^{i+k} \partial_i,
$
where $k=1,\ldots, n-2$, satisfy the relations
\begin{enumerate}
\item $\mathcal{L}_{X_{(k)}} g=0$ ~ (thus, they are isometries of $g$),
\item $\mathcal{L}_{X_{(k)}} \mu^{(n;\alpha)}=p_{[n,k,\alpha]} \mu^{(n;\alpha+k)}$,
\item $\mathcal{L}^m_{X_{(k)}} \mu^{(n;\alpha)}=\left(\prod_{s=0}^{m-1} (p_{[n,k,\alpha]}-2ks)\right)  \mu^{(n;\alpha+mk)}$,
\end{enumerate}
where the coefficients $p_{[n,k,\alpha]}$ are defined as $p_{[n,k,\alpha]}=3k+1-n-2\alpha.$
\end{lemma}
The proof of this lemma is a straightforward computation.

Consider now the general solution (\ref{Jor_sol}). Note that in the non-constant eigenvalue case,  $\xi_0\ne 0$, one can eliminate the constant term $\tilde g_0$ by a translation of variables $u^i$. Let $S_0$ be the resulting $n-1$ parameter family of solutions, 
\begin{equation}\label{Jor_sol_2}
S_0=\sum_{i=0}^{n-2} \xi_i \mu^{(n;i)},
\end{equation}
and let $\mathbf{L}_{[k]}$ be the Lie series
$$
\mathbf{L}_{[k]}=\exp(t_k\mathcal{L}_{X_{(k)}})=\sum_{s\ge0} \frac{t_k^s}{s!} \mathcal{L}^s_{X_{(k)}}, 
$$
where $X_{(k)}$ are as in Lemma \ref{lemmaLie}.
We point out that, when applied to $\mu^{(n;k)}$ for $n$ fixed, $\mathbf{L}_{[k]}$ consists of a finite number of terms: recall that  $\mu^{(n;i)}=0$ for  $i>n-2$.


\begin{lemma}\label{Norm1}
If $\xi_0\ne 0$, then it can be set equal to one.
\end{lemma}

\proof
Let us consider the scaling transformation
$$
v^i=\gamma^{\frac{n+1}{2}-i} u^i,
$$
where $\gamma \neq0$ is an arbitrary constant. It is easy to see that this  preserves the form of $g$. Direct calculation gives
$$
\mu^{(n;k)ij}(u)  \partial_{u^i} \partial_{u^j}=\gamma^{\frac{n-1}{2}+k} \mu^{(n;k)ij}(v) \partial_{v^i} \partial_{v^j}.
$$
Thus, setting $\gamma=\xi_0^{-\frac{2}{n-1}}$, we can reduce the coefficient of $\mu^{(n;0)}$ to $1$.
\endproof

To finish the proof of Theorem \ref{Jor_thm} we need  the following

\begin{proposition}\label{Norm2}
Suppose $\xi_0\neq0$. Then
\begin{enumerate}
\item if $n \not\equiv 1\!\!\!\!  \mod3$, there exists an orthogonal transformation which brings the $(n-1)$-parameter  solution $S_0$ to $ \mu^{(n;0)}$;
\item if $n \equiv 1 \!\!\!\! \mod3$, $n\neq4$, there exists an orthogonal transformation which brings $S_0$ to the one-parameter family $\mu^{(n;0)}+\kappa \mu^{(n;\frac{n-1}{3})}$,
\end{enumerate}
where $\kappa$ is an arbitrary constant.
\end{proposition}

\proof
By  Lemma  \ref{Norm1} we can consider the family $S_0$ in the form
$$
S_0=\mu^{(n;0)}+ \sum_{i=1}^{n-2} \kappa_i \mu^{(n;i)},
$$
where $\kappa_i$ are arbitrary constant coefficients.
Suppose $n \not\equiv 1\!\!\!\! \mod3$, then  the coefficients $p_{[k,n,0]}$ defined in Lemma \ref{lemmaLie} do not vanish. Let us apply $\mathbf{L}_{[1]}$ to $S_0$ and look at  the coefficient of $\mu^{(n;1)}$:
\begin{gather*}
\mathbf{L}_{[1]} S_0=S_0+t_1 \mathcal{L}_{X_{(1)}} S_0+\frac{t_1^2}{2} \mathcal{L}^2_{X_{(1)}} S_0+ \ldots +\frac{t_1^{n-2}}{(n-2)!}  \mathcal{L}^{n-2}_{X_{(1)}} S_0\\
= \mu^{(n;0)} + \left(\kappa_1+t_1 p_{[n,1,0]} \right)  \mu^{(n;1)}+ \ldots
\end{gather*}
We can always choose $t_1$ such that the coefficient of $\mu^{(n;1)}$ is zero. Let us call $S_1$ the resulting  $(n-2)$-parameter  family:
$$
S_1=\mathbf{L}_{[1]} S_0 |_{t_1=-\frac{\kappa_1}{p_{[n,1,0]}}} = \mu^{(n;0)}+\sum_{i=2}^{n-2} \tilde{\kappa}_i \mu^{(n;i)}.
$$
Applying $\mathbf{L}_{[2]}$  and looking at the coefficient of $\mu^{(n;2)}$ we obtain
$$
\mathbf{L}_{[2]} S_1= \mu^{(n;0)} + \left(\tilde{\kappa}_2+t_2 p_{[n,2,0]} \right)  \mu^{(n;2)}+ \ldots
$$
Again, we can choose $t_2$ such that the coefficient of $ \mu^{(n;2)}$ vanishes, and so on. Ultimately, we get
$$
\mathbf{L}_{[n-2]} \mathbf{L}_{[n-3]} \cdots \mathbf{L}_{[1]} S_0=\mu^{(n;0)},
$$
as required.

To prove the second part of the proposition, let us set $n=3m+1$. It is easy to see that $\mathcal{L}_{X_{(m)}} \mu^{(n;0)}=0$, since the coefficient $p_{[n,m,0]}$ vanishes. Note that in this case $p_{[n,k,0]}\neq0$ for $k\neq m$. For  fixed $m$, until $k=m-1$ we can apply the same procedure as above, obtaining
$$
S_{m-1}= \mu^{(n;0)}+\sum_{i=m}^{n-2} \tilde{\kappa}_i \mu^{(n;i)}.
$$
At this point, applying $\mathcal{L}_{X_{(m)}}$ to $S_{m-1}$, we cannot eliminate the coefficient of $\mu^{(n;m)}$, since $\mathcal{L}_{X_{(m)}} \mu^{(n;0)}=0$. However, applying $\mathcal{L}_{X_{(m+1)}}$ and looking at the coefficient of $\mu^{(n;m+1)}$,
$$
\mathbf{L}_{[m+1]} S_{m-1}=\mu^{(n;0)} +\tilde{\kappa}_m \mu^{(n;m)} + \left(\tilde{\kappa}_{m+1}+t_{m+1} p_{[n,m+1,0]} \right)  \mu^{(n;m+1)}+ \ldots ,
$$
we can eliminate it. Following  the same method, we arrive at the canonical form
$$
S_{can}= \mu^{(n;0)} +\tilde{\kappa}_m \mu^{(n;m)}.
$$

\endproof

The case $n=4$ is special. Indeed, the $4 \times 4 $ Mokhov metric does not correspond to the single Jordan block case. In the $n=4$ Jordan block case, normal forms are presented in Section \ref{sect_4}.
This concludes the proof of Theorem \ref{Jor_thm}. In the constant eigenvalue case,  $\xi_0=0$, we cannot eliminate $\tilde{g}_0$ by a shift.   Finally, Theorem \ref{Jor_thm_2} can be established by following the same procedure as above.


\subsection{Relation to Frobenius manifolds}\label{sect_frobenius}
In this Section we demonstrate the relation between  Mokhov's Hamiltonian operator and the trivial Frobenius manifold associated with
 the cohomology ring of  projective space. For this purpose let us briefly recall the definition of a Frobenius manifolds.

\begin{definition}
A \emph{Frobenius manifold} $(M,g,\circ,e,E)$
 is a manifold $M$ endowed with a (pseudo)-metric $g$, a product $\circ$ on the tangent spaces $T_u M$ and  a pair of vector fields $e$ and  $E$ such that
\begin{itemize}
\item the product $\circ$ is commutative and  associative:
$$c^i_{jk}=c^i_{kj},\quad c^i_{jl}c^l_{kh}=c^i_{kl}c^l_{jh}.$$
\item $g$ is flat,  invariant with respect to $\circ$,
$$g_{ik}c^k_{jl}=g_{jk}c^k_{il},$$
 and the associated Levi-Civita connection $\nabla$ is compatible with $\circ$:
\begin{equation}
\nabla_lc^i_{jk}=\nabla_jc^i_{lk}.
\end{equation}
This implies that there exists a function $F$, called the \emph{Frobenius potential}, such that, in flat coordinates for $g$, 
$$c_{ijk}=g_{il}c^l_{jk}=\d_i\d_j\d_k F.$$
\item The product $\circ$ has a unity $e$ which is flat: $\nabla e=0$.
\item The Euler vector field $E$ satisfies: $$\nabla\nabla E=0,\quad{\rm Lie}_E e=-e,\quad{\rm Lie}_E \circ=\circ,\quad{\rm Lie}_E g=(2-d)g,$$
for some constant $d$. The existence of the Euler vector field is related to the existence of a flat
 contravariant metric called the \emph{intersection form}. In local coordinates it is defined by the formula
$$\tilde{g}^{ij}=g^{il}c^j_{lk}E^k.$$ 
\end{itemize}
\end{definition}

If, in the flat coordinates for $g$, the functions $c^i_{jk}$ are constant, the Frobenius manifold is called \emph{trivial} \cite{D}. In this case, the Frobenius potential in a cubic polynomial,  $F=\frac{1}{6}c_{ijk}u^iu^ju^k$.
 We now define  the trivial Frobenius manifold associated with Mokhov's Hamiltonian operator.
  
\begin{theorem}
The metric
\begin{displaymath}
g_{ij} = \left\{ \begin{array}{ll}
1 & \textrm{if $i+j=n+1$}\\
0 & \textrm{otherwise}
\end{array},\right.
\end{displaymath}
the structure constants
\begin{displaymath}
c^{i}_{jk} =g^{il}c_{ljk}=\left\{ \begin{array}{ll}
1 & \textrm{if $l=2n+1-j-k=n+1-i$ that is $j+k-i=n$}\\
0 & \textrm{otherwise}
\end{array},\right.
\end{displaymath}
the unity $e=\frac{\d}{\d u^n}$, 
and the Euler vector field $E=\sum_{k=1}^n(3k-2n-1)u^k\frac{\d}{\d u^k}$ define a trivial 
 Frobenius manifold with $d=3$. Moreover, the intersection form, 
\begin{displaymath}
\tilde{g}^{ij}=g^{il}c^j_{lk}E^k=\left\{ \begin{array}{ll}
[3(i+j)-2n-4]u^{i+j-1} & \textrm{if $i+j-1\le n$}\\
0 & \textrm{otherwise}
\end{array},\right.
\end{displaymath}
coincides with the second metric of Mokhov's operator.
\end{theorem}

The proof is a straightforward computation. 

\endproof

To conclude this section we compare the Frobenius algebra underlying Mokhov's example with the
 Frobenius algebra structure on the full cohomology ring of projective space $H^*({\bf CP}^d)$. This can be defined  with respect to the natural basis 
$$e_1=1,\ e_2=\omega,\ \dots,\ e_{d+1}=\omega^d,$$
generated by  powers of the standard Kahler form  normalized as
$$\int_{{\bf CP}^d}\omega^d=1.$$
The contravariant components of the metric $g$ and the structure constants
 $c^i_{jk}$ are defined respectively by
\begin{displaymath}
g_{ij} = \left\{ \begin{array}{ll}
1 & \textrm{if $i+j=n+1$}\\
0 & \textrm{otherwise}
\end{array},\right.
\end{displaymath}
and 
$$e_j\wedge e_k=c^i_{jk}e_i=e_{j+k-1},$$
that is,
\begin{displaymath}
c^{i}_{jk} =\left\{ \begin{array}{ll}
1 & \textrm{if $j+k-i=1$}\\
0 & \textrm{otherwise}
\end{array},\right.
\end{displaymath}
Putting $i'=n+1-i$ we obtain the Frobenius algebra of Mokhov's example.


\section{Hamiltonian operators in higher dimensions}\label{sect_3D}
In this section we consider general $d$-dimensional $n$-component Hamiltonian operators of the form
\begin{equation}\label{Ham_N}
P^{ij}=\sum_{\alpha=1}^{d} \left(g^{ij\alpha}({\bf u}) \frac{d}{dx^\alpha}+b^{ij \alpha}_k({\bf u}) u^k_{x^\alpha}\right).
\end{equation}
Operator of this type is called non-degenerate if  a generic linear combination of the bivectors $g^{\alpha}$ is non-degenerate (without any loss of generality we will assume that each  $g^{\alpha}$ is non-degenerate: this can always be achieved by a suitable linear transformation of the independent variables $x^{\alpha}$). In the non-degenerate  case, Mokhov's conditions involving the obstruction tensor must be satisfied by each pair of bivectors  (see \cite{Mokhov1} for further details). This implies that each pair  $(g^{\beta},g^{\gamma})$  defines a 2D Hamiltonian operator.
For a generic $d$-dimensional Hamiltonian operator, the Mokhov conditions   can be reformulated as follows
\begin{theorem}
Suppose $g^{\alpha}$ are flat contravariant metrics. An operator of the form \eqref{Ham_N}
defines a $d$-dimensional Hamiltonian operator if and only if for all $\beta\neq\gamma$ the following conditions are fulfilled: 
\begin{enumerate}
\item Linearity of $g^{\beta}$ in the flat coordinates of $g^{\gamma}$.
\item Vanishing of the Nijenhuis torsion of the affinor $L^{(\beta \gamma)}=g^{\beta} (g^{\gamma})^{-1}$.
\item Killing condition for $g^{\gamma}$ with respect to $g^{\beta}$: 
$$\nabla^i g^{kj\gamma}+\nabla^k g^{ij\gamma}
+\nabla^j g^{ik\gamma}=0,$$
where $\nabla$ is the Levi-Civita connection of $g^{\beta}$.
\end{enumerate}
\end{theorem}
\begin{remark}
It is sufficient to require the flatness of only one of the metrics $g^{\alpha}$. Indeed, let us suppose that $g^{\alpha}$ is flat. Then, since the pair $(g^{\alpha},g^{\beta})$  defines a 2D Hamiltonian operator for all $\beta$, the linearity, Nijenhuis and Killing conditions imply the flatness of $g^{\beta}$ (see Theorem 2).
\end{remark}

It was demonstrated by Mokhov \cite{Mokhov2} that there exist no non-constant 3D Hamiltonian operators with one or two components.
Here we show that there exist only two non-trivial three-component  Hamiltonian operators in 3D, namely

\medskip

\noindent {\bf Theorem 5.} {\it Any non-degenerate three-component Hamiltonian operator in 3D, which is not transformable to constant coefficients,  can be brought to one of the two canonical forms:
$$
P=
\begin{pmatrix}
\partial_z & 0 & \partial_x \\
0 & \partial_x & 0 \\
\partial_x & 0 & 0
\end{pmatrix}
+
\begin{pmatrix}
-2u^2 \partial_y -u_y^2 & u^3 \partial_y +2 u_y^3 & 0 \\
u^3 \partial_y -u_y^3 & 0 & 0\\
0 & 0 & 0
\end{pmatrix},
$$
or
$$
P=
\begin{pmatrix}
0 & \partial_x & 0 \\
\partial_x & 0 & 0 \\
0 & 0 & \partial_z
\end{pmatrix}
+
\begin{pmatrix}
-2 u^1 \partial_y -u_y^1  & u^2\partial_y +2 u_y^2 & 0 \\
u^2 \partial_y -u_y^2& 0 & 0\\
0 & 0 & 0
\end{pmatrix},
$$
by a local change of the dependent variables $u^i$, and a  linear  change of the independent  variables $x,y,z.$
Note that the second operator is reducible. Thus, there exists a unique three-component irreducible operator in 3D.}

\medskip

\proof
Since we are interested in the non-constant case, we will consider 3D operators as deformation of 2D non-constant operators which have been classified already (Theorem 3). There exist only three such operators, defined by the following pairs of contravariant metrics:
\begin{equation}\label{2D_1}
g=\begin{pmatrix}
0 & 0 & 1 \\
0 & 1 & 0 \\
1 & 0 & 0
\end{pmatrix},
\quad
\tilde g=
\begin{pmatrix}
-2 u^1 & -\frac{1}{2}u^2 & u^3\\
-\frac{1}{2}u^2 & u^3 & 0\\
 u^3 & 0 & 0
\end{pmatrix},
\end{equation}
\begin{equation}\label{2D_2}
g=\begin{pmatrix}
0 & 0 & 1 \\
0 & 1 & 0 \\
1 & 0 & 0
\end{pmatrix},
\quad
\tilde g=
\begin{pmatrix}
-2u^2 &u^3 & \lambda \\
u^3 & \lambda & 0\\
\lambda & 0 & 0
\end{pmatrix},
\end{equation}
and (reducible case)
\begin{equation}\label{2D_3}
g=
\begin{pmatrix}
0 & 1 & 0 \\
1 & 0 & 0 \\
0 & 0 & 1
\end{pmatrix},
\quad
\tilde g=
\begin{pmatrix}
-2 u^1 & u^2 & 0 \\
u^2 & 0 & 0\\
0 & 0 & \lambda
\end{pmatrix}.
\end{equation}
Fixing one of the above pairs,  let us add a third contravariant metric $h$. Since we are in the flat coordinates of the first metric $g$, the bivector $h$ must be linear.  Since the pair $(g,h)$ satisfies the Killing condition, we can represent  $h$ as a sum of symmetrized tensor products of infinitesimal  isometries  of $g$. Assuming this, let us consider  the  above three cases separately.

\noindent
{\bf Case \eqref{2D_1}:}
Checking the Killing condition for the pair $(\tilde g,h)$ we obtain that $h$ must be a linear combination of $g$ and $\tilde g$. This means that  our operator is essentially two-dimensional.

\noindent
{\bf Case \eqref{2D_2}:}
Checking the Killing condition for the pair $(\tilde g,h)$ we obtain
$
h=c_1 g +c_2 \tilde g +h_0
$
where
$$
h_0=
\begin{pmatrix}
\nu & 0 & 0 \\
0 & 0 & 0 \\
0 & 0 & 0
\end{pmatrix}.
$$
One can verify that this ansatz for $h$  satisfies all other conditions. Finally, $c_1, c_2$ and $\lambda$ can be eliminated, and $\nu$ can be set equal to $1$ by a linear change of the independent variables $x,y,z$. This gives the first (irreducible) case of Theorem 5.

\noindent
{\bf Case \eqref{2D_3}:}
In this case, it is no longer sufficient to consider  the Killing condition alone: we also need  the linearity of $h$ with respect to $\tilde g$, that is $\nabla^2 h=0$, where $\nabla$ correponds to $\tilde g$. These conditions imply $h=c_1 g +c_2 \tilde g +h_0$, where
$$
h_0=
\begin{pmatrix}
0 & \mu & 0 \\
\mu & 0 & 0 \\
0 & 0 & \nu
\end{pmatrix}.
$$
One can verify that this ansatz for $h$ satisfy all other conditions.
By a linear transformation of $x,y,z$ we can eliminated $c_1, c_2$ and $\lambda$, transforming the operator  to the second (reducible) case of Theorem 5. 
\endproof

It follows from the proof of Theorem 5 that any non-degenerate three-component Hamiltonian operator in 4D is essentially 3D, or can be transformed to constant coefficient form.  We point out that there exists non-trivial examples of Hamiltonian operators in any dimension:

\begin{example}
The following expression provides an  example of non-constant irreducible $N$-component Hamiltonian operator in $N$ dimensions:
$$
P^{ij}=\eta^{ij}\frac{d}{dx^1}+g^{ij}\frac{d}{dx^2}+b^{ij}_k u^k_{x^2}+\sum_{m=1}^{N-2} h^{ijm} \frac{d}{dx^{m+2}},
$$
where
\begin{itemize}
\item $\eta$ is the constant $N \times N$ anti-diagonal bivector: $\eta^{ij}=\delta^{i,N+1-j}$;
\item the bivector $g$  is defined as $g^{ij}=\mu^{(N;N-2)}+g^{ij}_0$, where $\mu^{(N;N-2)}$ is defined by \eqref{Mok_shift}, namely
$$
\mu^{(N;N-2) ij}=[3(i+j)-8] u^{i+j+N-3},
$$
and $g^{ij}_0=\delta^{i,N-j}+\lambda \delta^{i,N+1-j}$;
\item $b^{ij}_k$ are the contravariant Christoffel symbols of $g$, namely
$$
b^{11}_{N-1}=-1, \quad b^{12}_{N}=2, \quad b^{21}_{N}=-1,
$$
with all  remaining coefficients equal to $0$;
\item the $N-2$ constant bivectors $h^m$ are defined as $h^{ijm}=\delta^{i m} \delta_{m}^{j}$.
\end{itemize}
For $N=3$, this example corresponds to the first case of Theorem 5.
\end{example}


\section{Concluding remarks}

This paper outlines an approach to the classification of  first order multi-dimensional Hamiltonian operators of differential-geometric type.  Our main contributions include a complete list of 2D Hamiltonian operators with three and four components, as well as the classification of  multi-component operators in case where the corresponding affinor consists of Jordan blocks with distinct eigenvalues.  

\begin{enumerate}

\item Our calculations demonstrate that the most challenging case is the one where the  Jordan normal form of the affinor $L=\tilde g g^{-1}$ consists of several Jordan blocks with the same eigenvalue. To complete the classification, one needs  to understand the structure of such operators: due to the splitting lemma, the general Hamiltonian operator would be representable as their direct sum.  

\item Given any 2D Hamiltonian operator from our list, it would be interesting to  classify Hamiltonians which generate  integrable $2+1$ dimensional systems of hydrodynamic type. The existing results suggest that `integrable' Hamiltonians form finite-dimensional moduli spaces, and are quite non-trivial even for constant-coefficient operators, see  \cite{Fer1, FMS, FOS} for the first steps in this direction. 

\item It would be interesting to develop a deformation theory of 2D Hamiltonian operators in the spirit of \cite{Getzler, Magri, DZ}, and to investigate  triviality of  Poisson cohomogy in 2D.
 Some results in this direction were recently obtained in \cite{C}. 
\end{enumerate}

\addcontentsline{toc}{section}{Acknowledgments}
\section*{Acknowledgments}
We  thank Alexey Bolsinov for numerous helpful discussions.
The research of PL was partially supported by the London Mathematical Society, by the Italian MIUR Research Project \emph{Teorie geometriche e analitiche dei sistemi Hamiltoniani in dimensioni finite e infinite}, and by  the GNFM Programme \emph{Progetto Giovani 2012}.



\end{document}